\documentclass[1p]{elsarticle}

\usepackage{amsmath,mathtools}
\usepackage{amssymb}
\usepackage{latexsym}
\usepackage{bm,bbm,mathrsfs}

\usepackage{tikz} 
\usepackage{pgfplots}
\usetikzlibrary{shapes,arrows,shadows,pgfplots.groupplots,decorations.markings,chains}
\pgfplotsset{compat=newest} 
\pgfplotsset{  tick label style={font=\scriptsize}, %
                            label style={font=\scriptsize}, %
                            legend style={nodes=right, font=\small}, %
                            title style={font=\scriptsize} ,%
                             xlabel style={yshift=5pt}} 

\DeclareMathOperator*{\argmin}{argmin}

\newcommand{\dsp}			{\displaystyle}
\newcommand{\R}				{\mathbbm{R}}
\newcommand{\1}				{\mathbbm{1}}

\newcommand{\M}				{\mathbbm M}

\newcommand{\vareps}		{\varepsilon}
\newcommand{\D}				{\scriptstyle{\mathrm{D}}}

\usepackage{amsthm}

\newtheoremstyle{Myremark}
{\topsep}
{\topsep}
{}
{}
{\bfseries \scshape}
{}
{\newline}
{\normalsize{\thmname{#1}\thmnumber{ #2}\thmnote{ (\scshape{#3})}}}

\newtheoremstyle{MyThm}
{\topsep}
{\topsep}
{\slshape}
{}
{\bfseries \scshape}
{}
{\newline}
{\thmname{#1}\thmnumber{ #2}\thmnote{ (\scshape{#3})}}

\theoremstyle{MyThm}

\theoremstyle{Myremark}
\newtheorem{Rem}{Remark}

\newenvironment{itemize*}%
  { \begin{itemize}
  \setlength{\itemsep}{1pt}
  \setlength{\parskip}{0pt}
  \setlength{\parsep}{0pt}}
  {\end{itemize}}

\usepackage[normalem]{ulem}
\definecolor{myorange}{rgb}{0.9568,0.4941,0.1961}
\definecolor{myred}{rgb}{0.9098,0.1294,0.2078}
\definecolor{myblue}{rgb}{0.0352,0.4981,0.6509}
\definecolor{mygreen}{rgb}{0.2235,0.6353,0.2588}
\definecolor{lightgray}{rgb}{0.8,0.8,0.8}

\usepackage{manfnt}

\newcommand{\highlight}[1]{{\color{myblue}#1}}

\newcommand{\mech}								{\mathrm{m}}

\newcommand{\Dist}		 					    {\mathrm{dist}}
\newcommand{\diff}[0]	 					    {\mathrm{d}}

\newcommand{\dpart}[2]							{\frac{\partial #1}{\partial #2}}

\newcommand{\Grad}		 					    {\vect{\nabla}\,}

\newcommand{\Div}								{\vect{\nabla} \cdot}

\newcommand{\abs}[1]	 					    {\left| #1 \right|}
\newcommand{\norm}[1]	 					    {\lVert#1\rVert}
\newcommand{\tens}[1]	 					    {\underline{\underline{#1}}}
\newcommand{\vect}[1]	 					    {\underline{#1}}

\newcommand{\domain}	 					    {\Omega}
\newcommand{\domainHeart}	 					{\domain^{\scriptscriptstyle \mathrm{H}}}
\newcommand{\x}			 					    {\vect{\mathrm{x}}}
\newcommand{\xRef}			 					{\vect{\xi}}

\newcommand{\domainRef}							{\domain_0}
\newcommand{\domainHeartRef}					{\domainHeart_0}

\newcommand{\disp}								{\vect{y}}
\newcommand{\vel}								{\vect{v}}

\newcommand{\velTest}							{\vect{v}^\natural}

\newcommand{\normal}							{\vect{n}}

\newcommand{\deforMap}							{\vect{\phi}}
\newcommand{\gradDef}						    {\tens{F}}
\newcommand{\stressCauchy}						{\tens{\sigma}}
\newcommand{\firstPiola}						{\tens{T}}
\newcommand{\secondPiola}						{\tens{\Sigma}}

\newcommand{\greenLagrange}						{\tens{e}}

\newcommand{\cauchyGreen}						{\tens{C}}

\newcommand{\strainLinear}						{\tens{\vareps}}

\newcommand{\elasticityTensor}					{\tens{\tens{A}}}

\newcommand{\hyperEnergy}						{\mathscr{W}}

\newcommand{\volmass}							{\rho}
\newcommand{\fiberStress}						{\sigma_{{\scriptscriptstyle 1\text{D}}}}
\newcommand{\fiber}								{\vect{\tau}}
\newcommand{\fiberStrain}						{\fiber \cdot \greenLagrange \cdot \fiber}
\newcommand{\activeStrain}						{e_c}
	
\newcommand{\Ndof}								{N}
\renewcommand{\P}								{\mathbb{P}}
\newcommand{\spaceDof}							{\mathcal{V}}
\newcommand{\dispDof}							{\vec{Y}}
\newcommand{\velDof}							{\vec{V}}

\newcommand{\velTestDof}						{\vec{V}^{\natural}}
\newcommand{\velTestDofintercal}				{\vec{V}^{\natural^{\scriptstyle \intercal}}}
\newcommand{\massDof}							{\bold{M}^{\rm m}}

\newcommand{\stressDof}							{\vec{\rm K}^{\rm m}}
\newcommand{\followPresureDof}					{\vec{\rm N}^{\rm m}}
\newcommand{\BestelDof}							{\vec{\rm B}^{\rm m}}

\newcommand{\domainMeasMech}					{\omega^\text{obs}_\mech}

\newcommand{\ext}								{\text{Ext}}
\newcommand{\extDof}							{\vec{\text{Ext}}}

\newcommand{\elec}								{\mathrm e}

\newcommand{\domainTorso}						{\Omega^{\scriptscriptstyle \mathrm{T}}}
\newcommand{\extCurrent}						{I_{\mathrm {app}}}
\newcommand{\extCurrentVec}						{\vec{\rm I}_{\mathrm {app}}}
\newcommand{\ionCurrent}						{I_{\mathrm{ion}}}
\newcommand{\ionCurrentVec}					{\vec{\rm I}_{\mathrm{ion}}}
\newcommand{\diffusion}   						{\tens{D}}
\newcommand{\intraDiffusion}   					{\diffusion_{\mathrm i}}
\newcommand{\extraDiffusion}   					{\diffusion_{\mathrm e}}
\newcommand{\torsoDiffusion}   					{\diffusion_{\mathrm T}}
\newcommand{\membArea}							{a_{m}}
\newcommand{\membCapacit}						{c_{m}}
\newcommand{\intraV}							{u_{\mathrm i}}
\newcommand{\extraV}							{u_{\mathrm e}}
\newcommand{\torsoV}							{u_{\mathrm T}}
\newcommand{\diffV}								{v_{m}}
\newcommand{\testV}								{u^\natural}

\newcommand{\internElec}						{w}
\newcommand{\vmax}   							{V_{\mathrm{max}}}
\newcommand{\vmin}  							{V_{\mathrm{min}}}
\newcommand{\vgate} 							{V_{\mathrm{gate}}}
\newcommand{\tin}   							{\tau_{\mathrm{in}}}
\newcommand{\tout}  							{\tau_{\mathrm{out}}}
\newcommand{\topen} 							{\tau_{\mathrm{open}}}
\newcommand{\tclose}							{\tau_{\mathrm{close}}}

\newcommand{\extraVDof}							{\vec{U}_{\mathrm e}}

\newcommand{\diffVDof}							{\vec{V}_{m}}
\newcommand{\dtdiffVDof}						{\dot{\vec{V}}_{m}}
\newcommand{\internElecDof}						{\vec{W}}
\newcommand{\testVDof}							{\vec{U}^\natural}

\newcommand{\scalMassDof}						{\bold{M}^{\elec}}
\newcommand{\laplaceDof}						{\bold{K}^{\elec}}

\newcommand{\reduced}							{{\scriptscriptstyle \mathrm{POD}}}
\newcommand{\basisPod}							{\Pi}
\newcommand{\dimPod}							{N_{\reduced}}
\newcommand{\massPOD}							{M_{\reduced}}

\newcommand{\state}			   					{x}
\newcommand{\stateDof}		   					{X}
\newcommand{\Aug}								{{}^{\flat\!}}

\newcommand{\Nstate}		   					{N_{\!{\scriptscriptstyle\stateDof}}}

\newcommand{\observ}		   					{z}
\newcommand{\observDof}		   					{Z}

\newcommand{\param}			   					{\vartheta}
\newcommand{\paramDof}		   					{\theta}
\newcommand{\Nparam}							{d}
\newcommand{\target}		   					{\dag}
\newcommand{\stateInit}							{\state_\diamond}
\newcommand{\stateInitDof}						{\stateDof_\diamond}
\newcommand{\paramInit}		   					{\param_\diamond}
\newcommand{\paramInitDof}		   				{\paramDof_\diamond}

\newcommand{\modelOp}		   					{A}

\newcommand{\obsOp}			   					{H}

\newcommand{\gainFilter}	   					{G}

\newcommand{\modelOpDof}	   					{\mathrm{\modelOp}}

\newcommand{\obsOpDof}		   					{\mathrm{\obsOp}}

\newcommand{\gainFilterDof}	   					{\mathrm{\gainFilter}}                              					
\newcommand{\stateGainFilterDof}				{\mathrm{\gainFilter}^{\!{\scriptscriptstyle\stateDof}}}
\newcommand{\paramGainFilterDof}				{\mathrm{\gainFilter}^{\paramDof}}
\newcommand{\errorStateLinDof}					{\delta\!\tilde{\stateDof}}

\newcommand{\obsNoise}		   					{\chi}

\newcommand{\initNoise}		   					{\zeta}
                               					
\newcommand{\obsNoiseNorm}	   					{M}

\newcommand{\obsNoiseCov}	   					{W}

\newcommand{\Cov}			   					{P}

\newcommand{\initNoiseCov}	   					{\Cov_{\diamond}}
\newcommand{\paramNoiseCov}	   					{\Cov_{*}}

\newcommand{\initNoiseStateDof}  				{\initNoise^{\scriptscriptstyle \stateDof}}                            					
\newcommand{\initNoiseParamDof}  				{\initNoise^{\paramDof}}                            					
\newcommand{\obsNoiseNormDof}  					{\mathrm{\obsNoiseNorm}}

\newcommand{\obsNoiseCovDof}   					{\mathrm{\obsNoiseCov}}
\newcommand{\initNoiseCovDof}  					{\mathrm{\Cov}_{\diamond}}

\newcommand{\CovDof}		   					{\mathrm{\Cov}}
\newcommand{\stateCovDof}		   				{\mathrm{\Cov}^{{\scriptscriptstyle\stateDof\! \stateDof}}}	
\newcommand{\stateCovDofPred}					{\mathrm{\Cov}^{{\scriptscriptstyle\stateDof\! \stateDof}-}}	
\newcommand{\stateCovDofCorr}					{\mathrm{\Cov}^{{\scriptscriptstyle\stateDof\! \stateDof}+}}
\newcommand{\paramCovDof}		   				{\mathrm{\Cov}^{\paramDof \paramDof}}											
\newcommand{\paramCovDofPred}					{\mathrm{\Cov}^{\paramDof \paramDof -}}	
\newcommand{\paramCovDofCorr}					{\mathrm{\Cov}^{\paramDof \paramDof +}}
\newcommand{\crossCovDof}		   				{\mathrm{\Cov}^{{\scriptscriptstyle \stateDof} \paramDof}}
\newcommand{\crossCovDofPred}					{\mathrm{\Cov}^{{\scriptscriptstyle\stateDof\!} \paramDof -}}	
\newcommand{\crossCovDofCorr}					{\mathrm{\Cov}^{{\scriptscriptstyle\stateDof\!} \paramDof +}}
\newcommand{\crossObsCovDof}		   			{\mathrm{\Cov}^{{\scriptscriptstyle \stateDof \! \observDof}}}
\newcommand{\crossObsCovParamDof}		   		{\mathrm{\Cov}^{{\scriptscriptstyle \paramDof \! \observDof}}}
\newcommand{\obsCovDof}		   					{\mathrm{\Cov}^{{\scriptscriptstyle \observDof \! \observDof}}}
\newcommand{\dotStateCovDof}		   			{\dot{\mathrm{\Cov}}^{{\scriptscriptstyle\stateDof\! \stateDof}}}	
\newcommand{\dotParamCovDof}		   			{\dot{\mathrm{\Cov}}^{\paramDof \paramDof}}											
\newcommand{\dotCrossCovDof}		   			{\dot{\mathrm{\Cov}}^{{\scriptscriptstyle \stateDof} \paramDof}}
\newcommand{\paramNoiseCovDof} 					{\mathrm{\Cov}_{*}}									
\newcommand{\diffStateDof}						{\diff_{\! \scriptscriptstyle{\stateDof}}}
\newcommand{\diffParamDof}						{\diff_{\paramDof}}

\newcommand{\dimReduce}                         {r}

\newcommand{\Nsigma}							{p}								 					
									
\newcommand{\redCovDof}							{\mathrm{U}}    
\newcommand{\projCovDof}						{\mathrm{L}}                                    
\newcommand{\projStateCovDof}		 	   		{\mathrm{L}^{\!\!{\scriptscriptstyle \stateDof}}}							
\newcommand{\projParamCovDof}	 	   			{\mathrm{L}^{\!\paramDof}}

\newcommand{\weightUKFmat}						{D}
\newcommand{\redCovDofSqrt}						{C}
\newcommand{\covObsUKF}							{\Gamma}

\newcommand{\projStateCovDofPerp}		 	   	{\mathrm{L}^{\!\!{\scriptscriptstyle \stateDof_{\perp}}}}		                               				
                               					
\newcommand{\projAlphaCovDof}		 	   		{\mathrm{L}^{\!\!{\scriptscriptstyle \alpha}}}		                               				                               					

\newlength{\figlength}
\newlength{\capslength}
\usepackage{color}
\usepackage{ifthen}

\title{Identification of weakly coupled multiphysics problems. Application to the inverse problem of electrocardiography}
\author{Cesare Corrado$^{1,2}$, Jean-Fr\'ed\'eric Gerbeau$^{1,2}$, Philippe Moireau${^{3*}}$,\\
\small{$^{1}$ Inria Paris-Rocquencourt, 78153 Le Chesnay, France} \\
\small{$^{2}$ Sorbonne Universit\'es UPMC Univ Paris 6, UMR 7598 Laboratoire Jacques-Louis Lions, 75005 Paris, France} \\ 
\small{$^{3}$ Inria Saclay Ile-de-France, 91120 Palaiseau, France}\\
\small{$^*$ Corresponding author}
}

\begin{document}

\begin{abstract}
This work addresses the inverse problem of electrocardiography from a new perspective, by combining electrical and mechanical measurements. Our strategy relies on the definition of a model of the electromechanical contraction which is registered on ECG data but also on measured mechanical displacements of the heart tissue typically extracted from medical images. In this respect, we establish in this work the convergence of a sequential estimator which combines for such coupled problems various state of the art sequential data assimilation methods in a unified consistent and efficient framework. Indeed we aggregate a Luenberger observer for the mechanical state and a Reduced Order Unscented Kalman Filter applied on the parameters to be identified and a POD projection of the electrical state. Then using synthetic data we show the benefits of our approach for the  estimation of the electrical state of the ventricles along the heart beat compared with more classical strategies which only consider an electrophysiological model with ECG measurements. Our numerical results actually show that the mechanical measurements improve the identifiability of the electrical problem allowing to reconstruct the electrical state of the coupled system more precisely. Therefore, this work is intended to be a first proof of concept, with theoretical justifications and numerical investigations, of the advantage of using available multi-modal observations for the estimation and identification of an electromechanical model of the heart.
\end{abstract}

\maketitle

\section{Introduction} 
\label{sec:introduction}

In the last few years, more and more attention has been paid to the problem of state and parameters identification for complex three-dimensional models used in biomedical applications. Several works can be cited:  for example in cardiac electrophysiology \cite{boulakia-schenone-gerbeau-12,burger-mardal-nielsen-10,chapelle:hal-00834397,nielsen-sundnes-mardal-06} or in cardiac mechanics \cite{chabiniok-moireau-ea-11,krishnarmurthy-mcculloch-kerckhoffs-13,Xi-ea-10}, or also in hemodynamics \cite{bertoglio-moireau-gerbeau-11,delia-mirabella-veneziani-12,Moireau-Bertoglio-2013,perego-veneziani-vergara-11}. The observations available in these contexts are often multiphysics since several modalities can be used simultaneously: electrocardiograms, electrograms, MRI, CT scan, flow measurements with ultrasound, pressure measurement with catheters, myocardium thickness measurement with piezoelectric sensors, etc. Up to now, the multiphysics nature of these problems has sometimes been taken into account in the direct models but it has rarely been used in the inverse problems.  

The first purpose of the present study is to exhibit an example where multiphysics observations actually improve the identifiability of a coupled system. More precisely, considering an electromechanical model of the heart, it is shown that the estimation of an electrical parameter is improved if the electrical observation, namely the electrocardiogram, is enriched with mechanical observations, namely the movement of the myocardium. In this model, the electrophysiology acts as an input for the mechanics, but the electromechanical feedback is neglected~\cite{chapelle:hal-00678772}. The coupling is therefore only one-way. The second purpose of this article is to show that a reduced order filtering strategy is well-suited to this class of multiphysics problems. Optimal filtering, like Kalman filter and its nonlinear extensions, is known to be very efficient but also too expensive to be used for large problems like those considered here. An effective strategy consists in using this kind of methods for the parameters only, and to address the uncertainties on the state variables through a less expensive approach. This strategy has been successfully used for example in mechanics where the state variable was handled with a Luenberger filter \cite{PM-DC-PLT-08}. In this paper, a similar strategy is proposed to filter the electrical state variables, by exploiting the special structure of the one-way coupled problem. But due to the structure of the equations, a Luenberger approach is not straightforward in electrophysiology. A different reduced order method, based on Proper Orthogonal Decomposition is then proposed, as it was done for cardiac mechanics in \cite{chapelle:hal-00834397}. 

\subsection{Background and related work}

The goal of the inverse problem of electrocardiography, also called {\em cardiac electrical imaging}, is to reconstruct the electrical activity of the heart from body surface potential maps. Various strategies have been proposed since four decades. All  of them assume that measurements of the electrical potential $\torsoV$ are available on parts of the torso boundary $\partial\domainTorso$. 

The different strategies can be distinguished by the cardiac electrical source models they rely on. One of the first approaches was to estimate equivalent electrical dipoles \cite{gulrajani-roberge-savard-84,mirvis-dowdie-zettergren-77}. Another popular approach is to estimate the heart surface potential, usually called epicardial potential (even though pericardial potential would be more appropriate as noted in \cite{messnarz-hanser-tilg-04}). The potential $\torsoV$ within the torso $\domainTorso$ is assumed to be solution of the Poisson problem:
\begin{equation}\label{eq:poisson}
	\begin{cases}
		\Div (\torsoDiffusion \cdot \Grad \torsoV) = 0, & \quad \text{in } \domainTorso\\ 
		\torsoV = \extraV, & \quad \text{on } \partial \domainHeartRef 
	\end{cases}
\end{equation}
where $\partial \domainHeartRef$ denotes the boundary of the heart and $\torsoDiffusion$ is the electrical conductivity of the torso. The inverse problem then consists in estimating  $\extraV$ on $\partial \domainHeartRef$   (see e.g. \cite{barr-ramsey-spach-77,burnes-taccardi-rudy-00,wang-kirby-johnson-11}). This problem being notoriously ill-posed, various regularizations have been proposed: Tikhonov \cite{messinger-rudy-88}, the use of temporal information \cite{greensite-huiskamp-98,oster-rudy-92}, truncated Singular Value Decomposition or truncated Total Least Square \cite{pullan-buist-cheng-05}.

Another approach to address the inverse problem of electrocardiography is to consider the following equation within the heart $\domainHeartRef$:
\[
-\Div ((\intraDiffusion+\extraDiffusion) \cdot \Grad \extraV) = \Div (\intraDiffusion \cdot \Grad \diffV),\quad \text{in } \domainHeartRef,
\]
which is one of two equations of the {\em bidomain model} (see Section~\ref{sub:electrophysiology}). Then, instead of estimating the epicardial potential $\extraV$ on the surface $\partial \domainHeartRef$, the goal is to estimate the transmembrane potential $\diffV$ within the heart. In \cite{messnarz-hanser-tilg-04}, the two approaches were compared: investigating the null-space of the inverse problem, the authors concluded that the transmembrane potential-based formulation is more promising because it is based on a stronger biophysical \emph{a priori}. As the epicardial potential approach, the transmembrane potential approach is ill-posed and must be regularized. In \cite{liu-liu-he-06}, four variants of $L^2$-Tikhonov regularization are compared. In \cite{nielsen-cai-lysaker-07}, $H^1$-Tikhonov regularization is used with a prior taking two different homogeneous values in the myocardium depending on cardiac phase (plateau or rest values). The inverse problem, formulated in a PDE-constrained framework, is  addressed by directly solving the optimality saddle-point problem. In \cite{nielsen-lysaker-tveito-07}, the estimation of the transmembrane potential is combined with a level set technique to efficiently identify the location of a myocardial infarction.
In \cite{wang-macleod-johnson-13}, the approach of \cite{nielsen-cai-lysaker-07} is generalized to more general objective functions and constraints in order to identify ischemic regions (characterized by lower amplitude during the plateau phase). Two different regularizations are investigated: the Tikhonov regularization, which is found to overestimate ischemic regions but with good sensitivity, and the Total Variation regularization which is found to underestimate ischemic regions but with high specificity.

In \cite{lopez-rincon-bendahmane-ainseba-13}, the authors note that the usual regularization techniques have no physical ground. Instead, they propose to regularize the inverse solution with the monodomain equations coupled to the Fenton-Karma ionic model. The strategy proposed in the present paper has some similarities with this approach. We rely on a full electrophysiological model of the action potential coupled to the Poisson problem \eqref{eq:poisson} to estimate the solution of the heart electrical activation. This physical model is personalized on the fly with respect to its parameters in order to adapt it to a specific patient. Furthermore, we propose an additional step of modeling by considering the mechanical response to the electrophysiological activation, so we are able to also integrate mechanical measurements. Indeed, we believe that multimodal observations improve the identifiability of the complete model and therefore improve the quality of the electrical and mechanical state reconstruction. Another originality of our work is the use of a sequential data assimilation strategy that is adapted to a coupled electromechanical evolution model. Here we demonstrate how state-of-the-art gain filter on the electrophysiological model and on the mechanical model can be aggregated to propose a joint gain filter for the coupled problem.

\subsection{Organization of the present work}

The paper is organized as follows. In Section~\ref{sec:models}, the electromechanical model is presented.  In Section~\ref{sec:measurements}, the observation operators -- namely the measurements -- are detailed for the electrical and the mechanical variables. In Section~\ref{sec:data_assimilation_principles}, the general notions of data assimilation, optimal filtering and reduced order filtering are reviewed. Although the algorithms of this section are not new, their presentation differs from what is most often done in the literature since a purely deterministic description is adopted. In Section~\ref{sec:data_assimilation_electromech}, the algorithms used for the electromechanical problem are proposed and analyzed. In Section~\ref{sec:numerical_test_cases}, numerical experiments based on synthetic data are presented. 
The main purpose is to estimate a non-homogeneous parameter of the electrical model using electrical and mechanical observations. 



\section{Models} 
\label{sec:models}

We present the models in a time and space continuous context before entering into the discretization details and their numerical implementation. Concerning the continuous context, we denote by an underline character any vector field of $\R^3$ and two underlines any second-order tensor.

\subsection{Electrophysiology} 
\label{sub:electrophysiology}

A widely accepted model of the macroscopic electrical activity of the heart is the so-called \emph{bidomain model}   \cite{CFS,sachse-04,sundnes-06,tung-78}. 
It consists of two degenerate parabolic reaction-diffusion PDEs which describe the dynamics of the averaged intra- and extracellular potentials $\intraV$ and $\extraV$, coupled to a system of ODEs defining an \emph{ionic model}. This model is related to the chemical dynamics of the myocardium cell membrane, in terms of the (vector or scalar) variable $\internElec$ representing the distributed ion concentrations and gating states, or a phenomenological counterpart. The model reads
\begin{equation}
	\label{eq:bidomain}
	\begin{cases}
		\membArea \Bigl( \membCapacit \partial_t \diffV + \ionCurrent(\diffV,\internElec) \Bigr) - \Div (\intraDiffusion  \Grad \intraV) = \membArea \extCurrent  ,&\quad\text{in } \domainHeartRef,\\
		\membArea \Bigl( \membCapacit \partial_t \diffV + \ionCurrent(\diffV,\internElec) \Bigr) + \Div (\extraDiffusion  \Grad \extraV) = \membArea \extCurrent ,&\quad \text{in } \domainHeartRef,\\
		\partial_t \internElec + g(\diffV,\internElec) = 0,&\quad\mbox{in } \domainHeartRef, 
	\end{cases}
\end{equation}
where $\diffV = \intraV - \extraV$ represents the transmembrane potential, $\membCapacit$ is the membrane capacitance per unit area, $\membArea$ is a constant representing the rate of membrane area per unit volume, $\intraDiffusion,\extraDiffusion$ are the intra- and extra-cellular conductivity tensors,  $\extCurrent$ is an external  volume current. In \eqref{eq:bidomain} the function $g$ represents an ionic model. In this article, the Mitchell-Schaeffer ionic model \cite{mitchell-schaeffer-03}  is considered, with the same rescaling as in \cite{boulakia:inria-00400490}. It is a reduced complexity model capable of integrating relevant phenomenological properties of the ventricle cell membrane: 
\begin{equation}\label{eq:mitchell}
		\begin{cases}
			\dsp \ionCurrent(\diffV,\internElec) = -\frac{\internElec}{\tin}\frac{(\diffV-\vmin)^2(\vmax-\diffV)}{\vmax-\vmin} +\frac{1}{\tout}\frac{\diffV-\vmin}{\vmax-\vmin}, \\[0.3cm]
			g(\diffV,\internElec) = \left\{
				\begin{aligned}
		      		\frac{\internElec}{\topen} -\frac{1}{\topen(\vmax-\vmin)^2}& \quad \mbox{if } \diffV \leq \vgate,\\
		      		\frac{\internElec}{\tclose} & \quad \mbox{if } \diffV > \vgate,
		  		\end{aligned}\right.
		\end{cases}
\end{equation}
where $\vgate$, $\tin$, $\tout$, $\topen$, $\tclose$ are given constants and $\vmin$ and $\vmax$ are scaling constants (typically $\vmin=-80 mV$ and $\vmax=20 mV$). 

On the boundary $\partial \domainHeartRef$, we have
$\intraDiffusion \cdot \Grad \intraV \cdot \normal = 0$, and the heart is assumed to be isolated, $\extraDiffusion \cdot \Grad \extraV \cdot \normal = 0$, as often done in the literature~\cite{clements-nenonen-li-horacek-04,potse-dube-gulrajani-06}. For well-posedness, the condition $\int_{\domainHeartRef} \extraV = 0$ is enforced. 

Summing and subtracting the first two equations of \eqref{eq:bidomain}, the system reads
\begin{equation}
	\label{eq:bidomain2}
	\begin{cases}
		\membArea \Bigl( \membCapacit \partial_t \diffV + \ionCurrent(\diffV,\internElec) \Bigr) \\
		\hspace{2cm}- \Div (\intraDiffusion \cdot \Grad (\extraV+\diffV)) = \membArea \extCurrent  ,&\quad\text{in } \domainHeartRef,\\
		-\Div ((\intraDiffusion+\extraDiffusion) \cdot \Grad \extraV) - \Div (\intraDiffusion \cdot \Grad \diffV) = 0 ,&\quad \text{in } \domainHeartRef,\\
		\partial_t \internElec + g(\diffV,\internElec) = 0, &\quad\mbox{in }  \domainHeartRef,  \\
		(\extraDiffusion \cdot \Grad \extraV) \cdot \normal = 0, & \quad \text{on } \partial \domainHeartRef, \\
		(\intraDiffusion \cdot \Grad (\extraV+\diffV)) \cdot \normal = 0, & \quad \text{on } \partial \domainHeartRef.
	\end{cases}
\end{equation}

A $\P_1$ finite element discretization of the potential variables is used --~see the corresponding refined electrical mesh in Figure~\ref{fig:meshes}. This leads to a discretization space $\spaceDof_h^\elec$ with $\Ndof^\elec$ degrees of freedom (i.e. $\spaceDof_h^\elec \simeq \R^{\Ndof^\elec}$). To each field --~for instance $\diffV$ --~is associated its corresponding approximation, denoted with an index $h$ --~for example ${\diffV}_h$. Equivalently, this approximation can be represented by its corresponding vector of degrees of freedom written with a vector in uppercase letter --~\emph{i.e.} $\diffVDof$. The finite element vector of degrees of freedom or the discretization of linear form are then defined with vectors in straight uppercase letter whereas the finite element matrix operator with a bold uppercase letter. With this notation: 
\[
	\forall \testVDof \in \spaceDof_h^\elec,\quad \testVDof  \scalMassDof \diffVDof = \int_{\domainRef} {\diffV}_h \testV_h \, d\Omega,\quad \testVDof  \scalMassDof_{a} \diffVDof = \int_{\domainRef} \membArea \membCapacit  {\diffV}_h \testV_h \, d\Omega,
\]
\[
	\forall \testVDof \in \spaceDof_h^\elec,\quad \testVDof \laplaceDof_{\rm i} \diffVDof = \int_{\domainRef} \intraDiffusion \cdot \Grad {\diffV}_h \cdot \Grad \testV_h \, d\Omega,
\]
\[
	\forall \testVDof \in \spaceDof_h^\elec,
	 \int_{\domainRef} \membArea \left(\extCurrent-\ionCurrent({\diffV}_h,\internElec_h)\right) \testV_h \, d\Omega
	\simeq  \testVDof \scalMassDof_{I}\left(\extCurrentVec - \ionCurrentVec(\diffVDof,\internElecDof)\right)
\]
and if a spatial discretization of the internal variable is done by node 
\[
	\forall \testVDof \in \spaceDof_h^\elec,\quad\int_{\domainRef} g({\diffV}_h,\internElec_h) \testV \, d\Omega
	  \simeq \testVDof \scalMassDof \vec{\rm G}(\diffVDof,\internElecDof)
\] 
where the applied current and the ionic variables have been interpolated, allowing to define the ionic variables at the nodes instead of the quadrature points.
Finally, after spatial discretization system \eqref{eq:bidomain2} reads 
\begin{equation}\label{eq:bidomain2Dof}
	\begin{cases}
\scalMassDof_a \dtdiffVDof + \laplaceDof_{\rm i} (\diffVDof + \extraVDof) = \scalMassDof_{I}(\extCurrentVec-\ionCurrentVec(\diffVDof,\internElecDof)),\\
		(\laplaceDof_{\rm i} + \laplaceDof_{\rm e}) \extraVDof + \laplaceDof_{\rm i} \diffVDof = 0, \\
		 \dot{\internElecDof} +\vec{\rm G}(\diffVDof,\internElecDof)= 0.
	\end{cases}
\end{equation}
The state of this system is $\stateDof^\elec = \Big(\begin{smallmatrix} \diffVDof \\ \internElecDof \end{smallmatrix}\Big)$ while $\extraVDof$ appears as an auxiliary variable verifying the static equilibrium \eqref{eq:bidomain2Dof}$_2$ with $\diffVDof$. This can be summarized by
\begin{equation}\label{eq:modelElec}
	\begin{cases}
		\dot{\stateDof}^\elec = \modelOpDof^\elec(\stateDof^\elec,\paramDof^\elec), \\
		\stateDof^\elec(0) = \stateInitDof^\elec + \initNoise_{\scriptscriptstyle \stateDof^\elec}, \\
		\paramDof^\elec = \paramInitDof^\elec + \initNoise_{\paramDof^\elec}.
	\end{cases}
\end{equation}
where  $\paramDof^\elec$ denotes the vector of parameters of the electrophysiology model, $\initNoise_{\paramDof^\elec}$ and $\initNoise_{\scriptscriptstyle \stateDof^\elec}$ the uncertainty on the parameters and the initial condition respectively.

\subsection{Mechanics} 
\label{sub:mechanics}

The heart domain is denoted by $\domainHeart(t)$ at any time $t$. This domain is the image of a reference configuration $\domainHeartRef$ through the solid deformation mapping $\deforMap$ 
\[
	\deforMap: 
	\left|
	\begin{aligned}
		\domainHeartRef \times [0,T] &\longrightarrow \domainHeart(t),\\
		(\xRef,t) &\longmapsto \x = \xRef + \disp (\xRef,t)	
	\end{aligned}
	\right.
\]
where $\disp$ is the solid displacement. The solid velocity is given by $\vel = \dot{\disp}$. The deformation gradient $\gradDef$ is given by
$\gradDef(\xi,t) = \tens{\nabla}_{\xRef} \vect{\phi} = \tens\1 + \tens{\nabla}_{\xRef}\disp$, and its determinant is denoted by $J$.
The right Cauchy-Green deformation tensor is defined by $\cauchyGreen = \gradDef^{\intercal}\cdot \gradDef$, the Green-Lagrange tensor by
$
	\greenLagrange = \frac{1}2 (\cauchyGreen-\tens\1) = \frac{1}2 \bigl(\tens{\nabla}_{\xRef}\disp + (\tens{\nabla}_{\xRef}\disp)^\intercal + (\tens{\nabla}_{\xRef}\disp)^\intercal \cdot \tens{\nabla}_{\xRef}\disp\bigr),
$
and its linearization by
$\strainLinear = \frac{1}2 \bigl(\tens{\nabla}_{\xRef}\disp + (\tens{\nabla}_{\xRef}\disp)^\intercal\bigr)$.

The mass per unit volume is denoted by  $\volmass$ and the Cauchy stress tensor by $\stressCauchy$. In the reference configuration, the first and second Piola stress tensor are respectively defined by $\firstPiola = J \stressCauchy \cdot \gradDef^{-\intercal}$ and $\secondPiola = \gradDef^{-1} \cdot \firstPiola = J \gradDef^{-1} \cdot \stressCauchy \cdot \gradDef^{-\intercal}$. The constitutive law is assumed to be a combination of a hyperelastic law of potential $\hyperEnergy$, a viscous component chosen proportional to the strain rate $\dot{\greenLagrange}$, and an active part along the fiber direction $\tau$ represented by 3 internal variables which are $e_c$ the active strain, $k_c$ the active stiffness and $\tau_c$ the associated active stress \cite{chapelle:hal-00678772}:
\begin{equation}\label{eq:constitutuvelaw}
		\secondPiola(\greenLagrange,e_c,k_c,\tau_c) = \dpart{\hyperEnergy}{\greenLagrange}(\greenLagrange) + \eta_s \dot{\greenLagrange} + \fiberStress(e_c,k_c,\tau_c) \fiber \otimes \fiber,
\end{equation}
with $\fiberStress = \frac{1+2 \activeStrain}{1+2 \fiberStrain} (\tau_c + \mu \dot{e}_c)$, 
where these 3 internal variables rely on a chemically-controlled constitutive law describing the myofibre mechanics \cite{Bestel01,chapelle:hal-00678772,Huxley57}:
\begin{align} \label{bestel}
&  \left\{
    \begin{aligned}
		\partial_t k_c = - (\abs{u}+\alpha \abs{\dot{e}_c})k_c + n_0 k_0 \abs{u}_+ & \quad\mbox{in }\domainHeartRef \\
		\partial_t \tau_c = - (\abs{u}+\alpha \abs{\dot{e}_c})\tau_c + \dot{e}_c k_c + n_0 \sigma_0 \abs{u}_+ & \quad\mbox{in }\domainHeartRef \\ 
    \end{aligned}
  \right.
\end{align}
with $\alpha$, $k_0$, $\sigma_0$ given parameters, $n_0$ a function of $e_c$ accounting for the Frank-Starling effect and $u$ directly related to the electrical activity of the heart by
\[
\highlight{u(t) = a\diffV(t) + b}
\]
where $a$ and $b$ are two scaling parameters.

Concerning the boundary conditions, the external organs are modeled by visco-elastic boundary conditions on a sub-part of the epicardium: $\firstPiola \cdot \normal = k_s \disp + c_s \vel \text{ on } \Gamma_n(t)$. A uniform pressure is enforced on the left and right endocardium: 
$\stressCauchy \cdot \normal_t = p_{v,i} \normal_t  \text{ on } \Gamma_{n,i}(t), i=\{1,2\}$.
In summary, the mechanical problem reads
\begin{equation} \label{eq:solid}
    \begin{cases}
		\partial_t \disp = \vel, & \quad\mbox{in }\domainHeartRef \\
		\volmass \partial_t \vel  - \Div(\firstPiola)  =  0,& \quad\mbox{in }\domainHeartRef \\
		\partial_t k_c = - (\abs{u}+\alpha \abs{\dot{e}_c})k_c + n_0 k_0 \abs{u}_+, & \quad\mbox{in }\domainHeartRef \\
		\partial_t \tau_c = - (\abs{u}+\alpha \abs{\dot{e}_c})\tau_c + \dot{e}_c k_c + n_0 \sigma_0 \abs{u}_+, & \quad\mbox{in }\domainHeartRef \\
		\firstPiola \cdot \normal = k_s \disp + c_s \vel, & \quad \text{ on } \Gamma_{n} \\
		\firstPiola \cdot \normal = J p_{v,i} \gradDef^{-\intercal} \cdot \normal, & \quad\text{ on }  \Gamma_{c,i} \\
		\firstPiola \cdot \normal = 0, & \quad \text{ on }  \partial \domainHeartRef  \backslash ((\cup_i \Gamma_{c,i}) \cup \Gamma_n)
    \end{cases}
\end{equation}
with the constitutive law \eqref{eq:constitutuvelaw}.

\begin{Rem}
Several improvements can be formulated on this model: other active or passive constitutive laws, more sophisticated boundary conditions for the tethering of the myocardium or hemodynamics, {\em etc}. The estimation strategy presented here can be adapted to any of these improvements since we only refer to state-space model description to formulate our data assimilation methods.
\end{Rem}

The system is discretized with a $\P_1$ finite element --~see the corresponding mechanical mesh in Figure~\ref{fig:meshes} --~with a 5\% compressibility acceptance in order to avoid any numerical locking. Concerning the fibre directions, we prescribed them on each point of the mesh with an elevation angle varying from -60 degrees to 60 degrees through the myocardium thickness. The discrete system is based on the variational formulation associated with \eqref{eq:solid}:
\begin{multline} \label{eq:solidVaria}
		\forall \velTest \in \mathcal{V}^v, \int_{\domainRef} \volmass \partial_t \vel \cdot \velTest \, d\Omega +  \int_{\domainHeartRef} \secondPiola(\greenLagrange,\dot{\greenLagrange},e_c,k_c,\tau_c) : \diff_{\disp} e \cdot \velTest \, d\Omega \\+  \int_{\Gamma_{n}} (k_s \disp + c_s \vel) \cdot \velTest \, d\Gamma = - \sum_i \int_{\Gamma_{c,i}} J p_v (\gradDef^{-\intercal} \cdot \normal) \cdot \velTest \, d\Gamma 
\end{multline}
Using the same convention as for the electrophysiological model discretization, the mass operator is defined by
\[
	\forall \velTestDof,  \velTestDofintercal \massDof \dot{\velDof} = \int_{\domainHeartRef} \volmass \partial_t \vel_h \cdot \velTest_h \, d\Omega,
\]
the stress residual by
\[
	\forall \velTestDof,\quad \velTestDofintercal \stressDof(\dispDof,\velDof) = \int_{\domainHeartRef} \secondPiola(\greenLagrange_h,\dot{\greenLagrange}_h,e_{c,h},k_{c,h},\tau_{c,h}) : \diff_{\disp} \greenLagrange_h \cdot \velTest_h \, d\Omega,
\]
some weighted mass operators on the boundary by
\[
	\forall \velTestDof,\quad \velTestDofintercal \massDof_{k_s,\Gamma_n} \dispDof = \int_{\Gamma_{n}} k_s \disp_h \cdot \velTest_h \, d\Gamma, \quad \velTestDofintercal \massDof_{c_s,\Gamma_n} \velDof = \int_{\Gamma_{n}} c_s \vel_h \cdot \velTest_h \, d\Gamma,
\]
and a following pressure operator by
\[
	\forall \velTestDof,\quad \velTestDofintercal \followPresureDof(\dispDof)  = \int_{\Gamma_{c,i}} J_h (\gradDef_h^{-\intercal} \cdot \normal) \cdot \velTest_h  \, d\Gamma .
\]
The internal variables $e_c,k_c,\tau_c$ are gathered in a vector $\imath_c$, which is discretized at the integration points. A vector of degrees of freedom $\vec{\imath}_c$ is associated with the discontinuous field  $(e_{c,h},k_{c,h},\tau_{c,h})$ and the model \eqref{bestel} is discretized by 
$\dot{\vec{\imath}}_c = \BestelDof(\vec{\imath}_c,\dispDof,\velDof)$. The complete mechanical model is spatially discretized into
\begin{equation} \label{eq:solidDof}
    \begin{cases}
		\dot{\dispDof} = \velDof,\\
		\massDof \dot{\velDof} + \stressDof(\dispDof,\velDof,\vec{\imath}_c) + \massDof_{k_s,\Gamma_n} \dispDof + \massDof_{c_s,\Gamma_n} \velDof =  - \sum_i p_{v,i} \followPresureDof(\dispDof),\\
		\dot{\vec{\imath}}_c = \BestelDof(\vec{\imath}_c,\dispDof,\velDof).
    \end{cases}
\end{equation}
The state of this system is $\stateDof^{\rm m} = \Bigl(\begin{smallmatrix} \dispDof \\ \velDof \\ \vec{\imath}_c \end{smallmatrix}\Bigr)$. 
Denoting by $\paramInitDof^{\rm m} $ the set of the parameters characterising the mechanics, affected by an a priori uncertainty $\initNoise_{\paramDof^{\rm m}}$ and with $\initNoise_{\scriptscriptstyle \stateDof^{\rm m}}$ the a priori uncertainty on the initial condition, the mechanical system reads
\begin{equation}\label{eq:modelMech}
	\begin{cases}
		\dot{\stateDof}^{\rm m} = \modelOpDof^{\rm m}(\stateDof^{\rm m},\stateDof^{\elec},\paramDof^{\rm m}), \\
		\stateDof^{\rm m}(0) = \stateInitDof^{\rm m} + \initNoise_{\scriptscriptstyle \stateDof^{\rm m}}, \\
		\paramDof^{\rm m} = \paramInitDof^{\rm m} + \initNoise_{\paramDof^{\rm m}},
	\end{cases}
\end{equation}
where the electrical variable $\stateDof^\elec$, solution to~\eqref{eq:bidomain2Dof}, can be seen as an input.


\subsection{Electromechanical coupling} 
\label{sec:application_to_a_coupled_electro_mechanical_problem}

From a computational point of view, the electromechanical problem consists in solving two coupled systems. We chose to keep the two sub-systems in independent solvers. This choice allows us to use legacy codes, to make their maintenance easier and to take advantage of the specific numerical methods adapted to each physical compartment. The coupling algorithm sketched in Figure~\ref{fig:masterworkers} is handled by a ``master'' code which exchanges the heart displacements and the transmembrane potential with the electrical and mechanical software. In this work, it is assumed that there is no electromechanical feedback. The transmembrane potential $\diffV$ is sent to the mechanical problem. Then, the one-way coupling is performed through the quantity $u = a\diffV + b$ which triggers the mechanical contraction via a change in the active stiffness $k_c$ and in the active stress $\tau_c$, see \eqref{bestel}.

\tikzstyle{block} = [rectangle, draw, text width=7em, text centered,drop shadow,fill=white, rounded corners, minimum height=3em]
\tikzstyle{pinstyle} = [pin edge={to-,thin,black}]
\pgfdeclarelayer{background}
\pgfdeclarelayer{foreground}
\pgfsetlayers{background,main,foreground}

\begin{figure}[htbp]
	\begin{center}
		
		\begin{tikzpicture}[auto, node distance=2cm,>=latex',scale=1]
			
		    \node [block] (MasterInit) {Initialization};		    
			\node [block, right of=MasterInit,node distance=4.45cm] (Mechanics0) {Mech. Solver \\ $t_n$};
			\node [block, left of=MasterInit, node distance=4.45cm] (Electric0) {Elec. Solver \\ $t_n$};
			\node [block, below of=MasterInit, node distance=4cm] (MasterIter) {Interpolation};
			
			\node [block, below of=Mechanics0,node distance=4cm] (Mechanics1) {Mech. Solver \\ $t_{n+1} = t_n + \Delta t_\mech$};
			\node [block,below of=Electric0, node distance=4cm] (Electric1) {Elec. Solver \\ $t_{n+1} = t_n + (\frac{\Delta t_\mech}{\Delta t_\elec}) \Delta t_\elec$};

			\draw [->,line width=1.1pt] (MasterIter.10) -- node[name=vToM] {$\diffV^{n+1}$} (Mechanics1.170);
			\draw [->,line width=1.1pt] (Electric1.10) -- node[name=vFromE] {$\diffV^{n+1}$}(MasterIter.170);
%
			
%
			\draw [->,line width=1.1pt] (MasterInit.350) -- node[name=dFromM] {\small Start Iter.\normalsize} (Mechanics0.190);			
			\draw [->,line width=1.1pt] (MasterInit.170) -- node[name=vFromE] {\small Start Iter. \normalsize} (Electric0.10);
			
			\node [below of=Electric0,node distance=1.4cm] (timeStep1) {$\bullet$};			
			\node [below of=Electric0,node distance=2cm] (timeStep2) {$\bullet$};
			\node [below of=Electric0,node distance=2.6cm] (timeStep3) {$\bullet$};
			\draw [->,line width=1.1pt] (Electric0.south) -- node[name=timeElec,above right] {$\Delta t_\elec$} (Electric1.north);
			
			\draw [->,line width=1.1pt] (Mechanics0.south) -- node[name=timeMech] {$\Delta t_\mech$} (Mechanics1.north);
			
			\draw [<-,line width=1.1pt] (MasterInit.south) -- node[name=update] {$n = n+1$} (MasterIter.north);
			
			\begin{pgfonlayer}{background}
			   \path (MasterInit.west)+(-0.2,2.1) node (MasterNW) {};
		       \path (MasterIter.east)+(0.2,-1.0) node (MasterSE) {};
			   \path[fill=black!20,rounded corners, draw=black!50, dashed]
			       (MasterNW) rectangle (MasterSE);
			   \path (MasterInit.north)+(0,1) node (MasterN) {Master Elec.-Mech.};
				
			    \path (Electric0.west)+(-0.2,2.1) node (ElectricNW) {};
		        \path (Electric1.east)+(0.2,-1.0) node (ElectricSE) {};

		        \path[fill=black!20,rounded corners, draw=black!50, dashed]
		            (ElectricNW) rectangle (ElectricSE);      
				\path (Electric0.north)+(0,1) node (ElectricN) {Slave Elec. Solver};

			    \path (Mechanics0.west)+(-0.2,2.1) node (MechanicsNW) {};
		        \path (Mechanics1.east)+(0.2,-1.0) node (MechanicsSE) {};

		        \path[fill=black!20,rounded corners, draw=black!50, dashed]
		            (MechanicsNW) rectangle (MechanicsSE);      
				\path (Mechanics0.north)+(0,1) node (MechanicsN) {Slave Mech. Solver};			
			 \end{pgfonlayer}				
		\end{tikzpicture}
		\caption{Master - Slave coupling} 
		\label{fig:masterworkers}
	\end{center}
\end{figure}
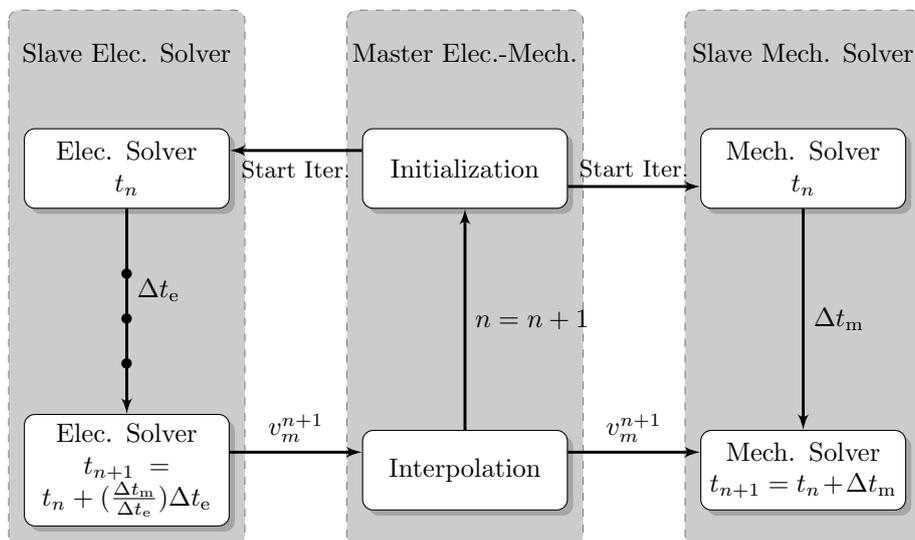

\begin{figure}
	[htb] 
	\begin{center} 
		\includegraphics[width=\textwidth]{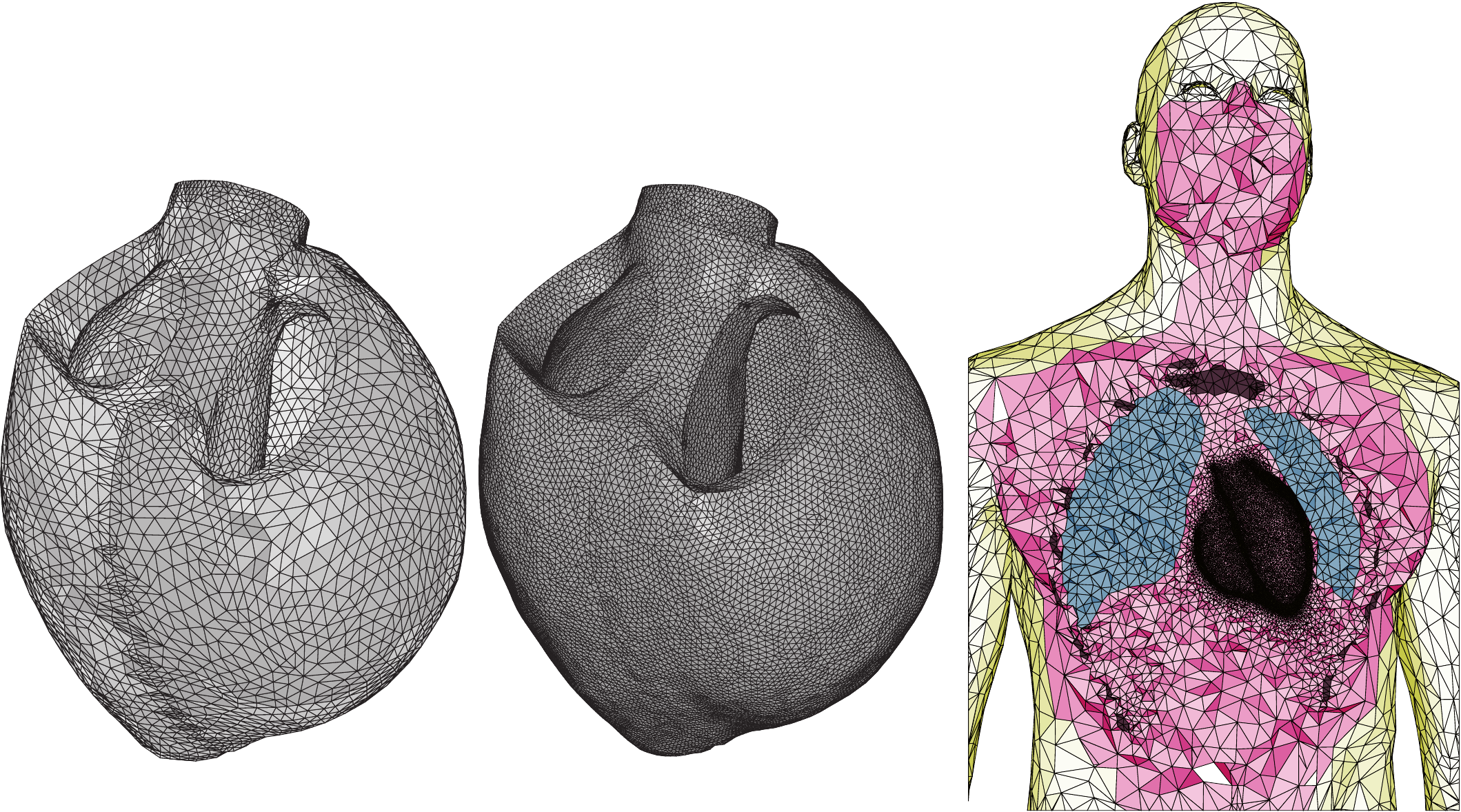}		
		 \caption{Mechanical mesh (Left, $4907$ nodes, $18193$ elements) and electrical mesh (center, $108112$ nodes, $541994$ elements) and thorax mesh (right, $229782$ nodes, $1250072$ elements)} \label{fig:meshes} 
	\end{center}
\end{figure}

\begin{figure}
	[htbp] 
	\begin{center} 
		\includegraphics[width=\textwidth]{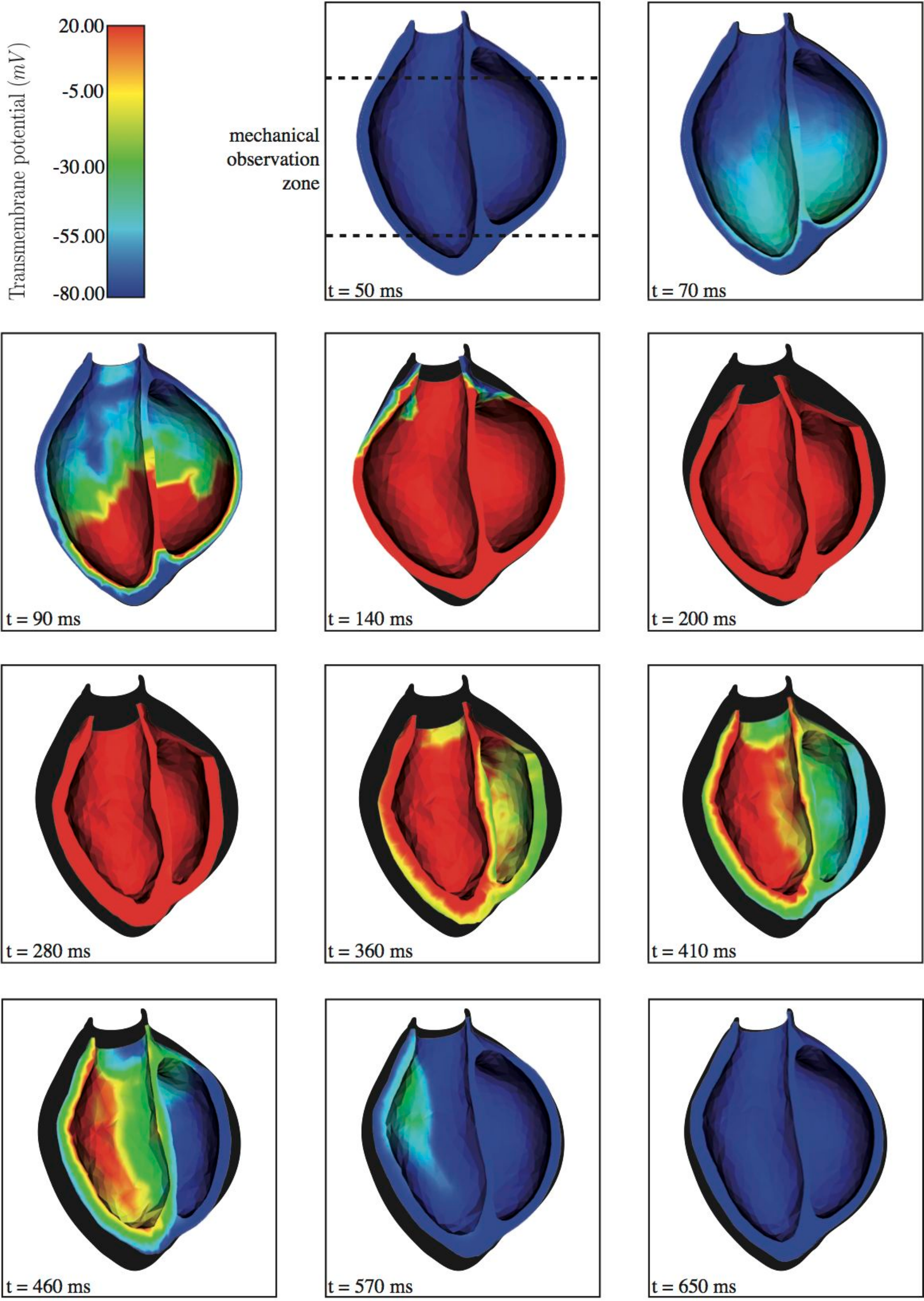}
		 \caption{Direct simulation results in long axis view --~deformed mesh with interpolated transmembrane potential from electrical refined mesh.} \label{fig:directsimulation} 
	\end{center}
\end{figure}

Because of the characteristics of the action potential propagation, the space and time steps required by the electrophysiology are typically smaller than those needed by the mechanics. To ensure accuracy at a reasonable computational cost, each sub-problem is solved with its own space and time step. The data are transmitted by the master code at some check-point in time (see Figure~\ref{fig:masterworkers}) and interpolated in space from a mesh to another (see Figure~\ref{fig:meshes}). Eventually, a typical complete direct simulation of our model is presented in Figure~\ref{fig:directsimulation}.

The inverse problem will be addressed on the discrete formulations. The state and the parameters are defined as the combination of the electrical and the mechanical ones
$
	\stateDof = \left(\begin{smallmatrix}
		\stateDof^\elec \\
		\stateDof^\mech
	\end{smallmatrix}\right), \quad 	
	\paramDof = \left(\begin{smallmatrix}
		\paramDof^\elec \\
		\paramDof^\mech
	\end{smallmatrix}\right).
$
Defining $\modelOpDof(\stateDof,\paramDof) = \left(\begin{smallmatrix}
		\modelOpDof^\elec(\stateDof^\elec,\paramDof^\elec) \\
		\modelOpDof^\mech(\stateDof^\mech,\stateDof^\elec,\paramDof^\mech)
	\end{smallmatrix}\right)$, 
this state variable follows the dynamics

\begin{equation}
\label{eq:electromech}
\dot{\stateDof} = \modelOpDof(\stateDof,\paramDof)
\end{equation}

The initial condition includes the quantities $\initNoiseStateDof$ and $\initNoiseParamDof$ which model the uncertainties on the initial condition and on the parameters respectively. When necessary, a trajectory will be denoted by $\stateDof_{[\initNoiseStateDof,\initNoiseParamDof]}$. This system can also be written in an augmented form by concatenated the state and parameters in $\Aug\stateDof = \bigl(\begin{smallmatrix} \stateDof \\ \paramDof \end{smallmatrix}\bigr)$. This augmented state satisfies $\Aug\dot{\stateDof} = \Aug\modelOpDof(\stateDof)$. The real trajectory is assumed to be a solution of the model for a given value of the uncertainties and is given by $\Aug\stateDof^\target = \bigl(\begin{smallmatrix} \stateDof^\target  \\ \paramDof^\target  \end{smallmatrix}\bigr)$. 

\section{Measurements}
\label{sec:measurements}

\subsection{Electrophysiological measurements} 
\label{sub:electrophysiological_measurements}

The best way to access the electrophysiological quantities would be to directly measure the electrophysiological potential at the heart surface. This is possible with an invasive procedure where a basket of electrodes is immersed in the ventricle and records the potential close to the endocardium surface (\cite{eldar1997transcutaneous} e.g.). In that case, the electrophysiological measurements can be represented by
\begin{equation}\label{eq:obs_elec_cont}
	\observDof^\elec(t) = \obsOp^\elec \state^\elec + \obsNoise^\elec	
\end{equation}
where the observation operator $\obsOp^\elec$ interpolates the values at the basket electrodes and $\chi^\elec$ gathers all the measurement errors. After spatial discretization
\begin{equation}\label{eq:obs_elec_disc}
	\observDof^\elec(t) = \obsOpDof^\elec \stateDof^\elec + \obsNoise^\elec	
\end{equation}
where now $\chi^\elec$ takes also into account the finite element discretization error.

To avoid such invasive measurement procedures, it is preferable to rely on electrocardiograms. Physiologically, the electrical potential diffuses from the heart to the rest of the body through the pericardium. The ECG depicts the time course of standardized potential differences on the body surface. We refer to \cite{boulakia:inria-00400490} for examples of healthy and pathological ECG obtained with the electrophysiology model used in the present paper.

%

The electrical diffusion within the body $\domainTorso$ can be modelled by a Poisson equation with the extracellular potential as a boundary condition on the heart surface and homogeneous Neumann boundary elsewhere:
\begin{equation}\label{eq:torso}
	\begin{cases}
		\Div (\torsoDiffusion \cdot \Grad \torsoV) = 0, & \quad \text{in } \domainTorso\\ 
		\torsoV = \extraV, & \quad \text{on } \partial \domainHeartRef \\
		(\torsoDiffusion \cdot \Grad \torsoV) \cdot \normal = 0, & \quad  \text{on } \partial \domainTorso \backslash \partial \domainHeartRef\
	\end{cases}
\end{equation}
where the  conductivity tensor is assumed to be isotropic, but non-homogeneous to account for the different conductivities of the lungs, the bones and the rest of the body (see \cite{boulakia:inria-00400490}).

After spatial discretization --~see the corresponding thorax mesh in Figure~\ref{fig:meshes} --~the ECG measurements denoted by $\observDof^\elec$ are still related to the electrical state variable by an observation operator $\obsOp^\elec$ which is the composition of the discrete diffusion operator and a linear combination of the potential at 9 points of the body surface. Assuming that the model is accurate enough, the actual ECG can still be represented by \eqref{eq:obs_elec_cont} where $\chi^\elec$ gathers all the possible measurement errors, including the errors resulting from the modelling of the measurement procedure. After spatial discretization, a relation of the form \eqref{eq:obs_elec_disc} still holds, but with $\chi^\elec$ also accounting for the discretization error.



\subsection{Mechanical measurements} 
\label{sub:mechanical_measurements}

The heart contraction can also be perceived from a kinematic and even a mechanical perspective. In fact, the displacements of the cavity can be observed using imaging modalities, for example Cine-MRI sequence, CT-sequence or even Tagged-MRI sequence. Considering for example this last type of measurement and assuming adequate registration strategies allowing to extract $3\D$ displacements of material points \cite{SPAMM_MotEst_AKR_SR_2008}, we can assume that displacements are measured in a large part of the ventricles, typically between two short axis planes, the first one slightly lower than the base and the second one slightly higher than the apex --~see Figure~\ref{fig:directsimulation}.  With more standard sequences, like Cine-MRI or CT sequences, we can benefit from optical flow strategies to extract displacements most of the time on part of the boundary. In both cases, after discretization, we can assume the existence of an observation operator $\obsOpDof^\mech$ such that the measured displacements are given by
\[
	\observDof^\mech(t) = \obsOpDof^\mech \stateDof^\mech + \obsNoise^\mech,
\]
where $\chi^\mech$ is now a vector taking into account the discretization error.

\begin{Rem}[Generalized observation operators]
Extracting displacements using registration techniques is still a challenge and it often leads to the measurement of an apparent motion instead of the real motion. As a consequence, it can be more adequate to consider the structure measured in itself more than its collection of material points. This was done typically in \cite{PM-DC-PLT-09} where in Cine-MRI and CT sequences the contours of interest are extracted and compared to the contours produced by the simulation by the use of a discrepancy field of the form $D(\xRef,t) = \vect{\Dist}_{\mathcal{S}_t}(\x(\xRef,t),t)=\obsNoise(\xi,t)$. Hence, the closer to 0 the discrepancy is the more information we have extracted from the images. As demonstrated in \cite{PM-DC-PLT-09}, this type of discrepancy operator is a simple non-linear generalization of the discrepancy computed by subtracting to the actual measurements the one produced by the use of the observation operator defined below. 
This generalization has demonstrated its efficiency in several real cases of data assimilation investigations \cite{chabiniok-moireau-ea-11,Moireau-Bertoglio-2013}. But, for the sake of simplicity, and since this work is primarily aimed at providing a proof of concept, we limit ourselves to more simple observation operators. Note that the use of a non-linear observation operator is, in the end, always analyzed in the light of the corresponding linearized operator \cite{PM-DC-PLT-09}.
\end{Rem}

The stresses are more intricate to measure. For now, the most common measurement available is an intracavity pressure obtained after invasive catheterization or reconstructed from arterial blood pressure measurements. In this article we assume that we can consider the intracavity pressure as a given source term in our model \eqref{eq:solid}. In a more general case, we could model the evolution of the intracavity pressures with the help for example of a Windkessel model and then use the possibly noisy measurements in a data assimilation context where the intracavity pressures are additional variables of the model.

\subsection{Multi-modalities measurements}
Finally the observations are concatenated and therefore given and the complete form
\begin{equation}\label{eq:obsOpTotal}
	\observDof = \begin{pmatrix}
		\observDof^\elec \\
		\observDof^\mech
	\end{pmatrix} = 	
	\begin{pmatrix}
		\obsOpDof^\elec & 0 \\
		0 & \obsOpDof^\mech
	\end{pmatrix}	
	\begin{pmatrix}
		\stateDof^{\elec\target} \\
		\stateDof^{\mech\target}
	\end{pmatrix} + 
	\begin{pmatrix}
		\obsNoise^\elec \\
		\obsNoise^\mech
	\end{pmatrix}
	= \obsOpDof \stateDof^\target + \obsNoise,
\end{equation}
which for the sake of generality will not be necessarily considered as linear when it is not mandatory. 


\section{Data assimilation principles}
\label{sec:data_assimilation_principles}

Data assimilation has become a very popular strategy to estimate a wide range of modeling uncertainties in numerical simulations \cite{Blum2008}. It has been initiated for environmental sciences but has now reached life sciences and especially cardiology. To introduce its main concepts in the \emph{state-space} formalism,  we consider a dynamical model
\begin{equation}\label{eq:modelNotation}
	\dot{\state} = \modelOp(\state,\param,t).
\end{equation}
In this equation, $\state$ denotes the \emph{state variable}, namely, the physical quantity which the model aims at describing during its time evolution. In this generic notation, the whole model is essentially summarized in the \emph{dynamical operator} $\modelOp$, which applies on the state variable itself, and may depend on time $t$ as well as on a set of physical parameters denoted by $\param$. In this work, \eqref{eq:modelNotation} corresponds to the electromechanical system defined by \eqref{eq:bidomain2} and \eqref{eq:solid}, or -- to avoid  issues associated with PDEs -- its discrete counterpart \eqref{eq:modelElec} and \eqref{eq:modelMech} (but most of what is presented here can be actually generalized to infinite dimensional systems~\cite{Bensoussan1971}). 

The initial condition $\state(0)$ and the parameter vector $\param$ must be prescribed in \eqref{eq:modelNotation}. Hence, when these quantities are unknown, they have to be estimated. In general $\state(0)$ and $\param$ are decomposed into known parts $\stateInit$ and $\paramInit$ called \emph{a priori} and uncertain parts $\initNoise^\state$ and $\initNoise^\param$:
\[
	\begin{cases}
		\state(0) = \stateInit + \initNoise^\state,\\
		\param = \paramInit + \initNoise^\param.
	\end{cases}
\]
When necessary, the trajectory  of $x$  will be denoted by $\state_{[\initNoise^\state,\initNoise^\param]}$.

Let us now re-consider an augmented state $\Aug\state = \left(\begin{smallmatrix} \state \\ \param \end{smallmatrix}\right)$ in order to integrate parameter uncertainties. The corresponding augmented model is defined by
\begin{equation}\label{eq:modelAugNotation}
	\Aug\dot{\state} = \Aug\modelOp(\Aug\state), \quad \text{ and } \Aug\state(0) = \begin{pmatrix} \stateInit + \initNoise^\state \\ \paramInit + \initNoise^\param \end{pmatrix} = \Aug\stateInit + \Aug\initNoise.
\end{equation}
Hence the parameters and the state can be consider in the same formalism. However, these two quantities differ in their dimension. Indeed, when the state is defined based on a PDEs model, the dimension of the state variable after space discretization is typically of the size of $10^3$ to $10^7$ degrees of freedom. By contrast, the size of the parameter vector is generally much more limited. Even if distributed parameters are considered, their variation should be considered smooth enough so that they can be discretized on a coarse mesh or a subdivision of the domain into large regions. Therefore, typically less than a hundred of parameter values need to be estimated. 

The measurements presented in Section~\ref{sec:measurements} can be cast in a general form: given a real trajectory $\state^\target$, the noisy observations are  represented using an observation operator $\obsOp$ such that
\begin{equation}\label{eq:measNotation}
	\observ = \obsOp(\state^\target,t) + \obsNoise(t),
\end{equation}
where $\observ$ denotes the actual data field. Note that $\obsOp$ can possibly be extended to a function of $\Aug\state$ and will then be denoted by $\Aug\obsOp$.

\subsection{Data assimilation by filtering} 
\label{sub:data_assimilation_by_filtering}

A common strategy in data assimilation is to minimize a criterion that balances the confidence in the a priori value of the state and parameters and a discrepancy measure between the given observations and the simulated ones. In the state-space formalism this criterion is typically least-square and reads
\begin{equation}\label{eq:criterion}
	\mathscr{J}(\initNoise) = \frac{1}2 \norm{{}^\flat \initNoise}_{\initNoiseCov^{-1}}^2 +\frac{1}2 \int_0^T \norm{\observ - \Aug\obsOp(\Aug\state_{[{}^\flat \initNoise]})}_\obsNoiseNorm^2 \, dt.
\end{equation}
where $\obsNoiseNorm$ is a metric on the observation space and $\initNoiseCov$ is the inverse of a metric which can be interpreted as an initial uncertainty covariance. The so-called \emph{4D-Var} method \cite{dimet2010variational} consists in finding the unknown quantities ${}^\flat\initNoise$ by minimizing $\mathscr{J}$ under the constraint of following the dynamics \eqref{eq:modelNotation}. When $T$ goes to infinity the minimization is expected to produce a better and better estimate. This minimization under constraint requires the computation of an adjoint variable used to compute the gradient. Hence, a gradient based descent algorithm requires numerous iterations of the direct and adjoint dynamics. 

When we do not wish to precisely retrieve the initial error ${}^\flat\initNoise$ but only seek to accurately approximate the state $\Aug\state=(\state,\param)$ ``independently'' of the possible initial error, we can avoid  minimization iterations by the use of sequential estimation methods. The principle in sequential estimation methods is to introduce a modified system --~denoted by a circumflex accent --~called \emph{observer} whose dynamics is changed to incorporate a correction based on the measured discrepancy. The new system dynamics therefore reads
\begin{equation}\label{eq:modelAugObserver}
	 \Aug\dot{\hat{\state}} = \modelOp(\Aug\hat{\state},t) + \gainFilter(\observ - \Aug\obsOp(\Aug\hat{\state})), \quad \Aug\hat{\state}(0) = \Aug\stateInit
\end{equation}
where $\gainFilterDof$ is called the \emph{observer gain} or \emph{filter}. 
The ultimate objective of the observer $\Aug\hat{\state}(t)$ is to converge in time to the real trajectory
\[
	\Aug\hat{\state}(t) \xrightarrow{t \to \infty} \Aug\state(t)
\]
and the gain has to be designed regarding this objective. Two strategies are commonly followed in this respect. First, the gain can be constructed from the optimality criterion \eqref{eq:criterion} by defining the so-called \emph{optimal observer}\footnote{the observer is called optimal due to the fact that it is associated with an optimal control problem.} --~or optimal sequential estimator --~with
\[
	\Aug\hat{\state}(t) = \Aug\state_{[\argmin \mathscr{J}(\cdot,t)]}(t).
\]
An optimal gain is obtained by differentiating this definition with respect to the time variable $t$ appearing in both the trajectory $\Aug\state(t)$ and in the criterion $\mathscr{J}(\cdot,t)$. This approach is well-known in a linear framework --~\emph{i.e.} where all the operators are linear --~and leads to the famous Kalman-Bucy filter \cite{Bensoussan1971,KalmanBucy61,Simon06}. The optimal filter is then defined from the solution of a Riccati equation. In a non-linear framework, the gain is more intricate to compute and derives from a Hamilton-Jacobi-Bellman solution \cite{Fleming:1997p148}. Hence, in this case, numerous works rely on approximate solutions as defined for example by the Extended Kalman Filter (EKF) which uses the Riccati equation of the linear Kalman filter with the tangent operator of the non-linear model and observation operators \cite{Simon06}. Eventually, the great advantage of the optimal sequential approach is that, like the 4D-Var, it can be defined for every model and every observation operator. The biggest drawback is that --~even with approximate solutions --~the filter is very costly to compute, especially for large dimensional systems like the ones produced by PDEs or their discretizations.  Therefore, alternative strategies -- the Luenberger observers and the Reduced Order Optimal filtering -- will be presented below.

\subsection{The Extended Kalman Filter (EKF)} 
\label{sub:the_extended_kalman_filter}

Even if we will have to rely on alternative strategies, it is interesting to develop the optimal filtering equations. For the sake of simplicity, we proceed after space discretization. The space discretized system variable is denoted by $\stateDof$ and the space discretized parameters by $\paramDof$. The observer state $\hat{\stateDof}$ and parameters $\hat{\paramDof}$ follow a modified version of the dynamics:
\begin{equation}\label{eq:observerGeneral}
	\begin{cases}
		\dot{\hat{\stateDof}} = \modelOpDof(\hat{\stateDof},\hat{\paramDof}) + \stateGainFilterDof (\observDof - \obsOpDof(\hat{\stateDof},t)), \\ 
		\dot{\hat{\paramDof}} = \paramGainFilterDof (\observDof - \obsOpDof(\hat{\stateDof},t)), \\
		\stateDof(0) = \stateInitDof, \\
		\paramDof(0) = \paramInitDof,
	\end{cases}
\end{equation}
where the gain $\gainFilterDof = \bigl(\begin{smallmatrix}
	\stateGainFilterDof \\ 
	\paramGainFilterDof
\end{smallmatrix}\bigr)$ is a linear operator. In the augmented form with $\Aug\hat{\stateDof} = \bigl(\begin{smallmatrix} \hat{\stateDof} \\ \hat{\paramDof} \end{smallmatrix}\bigr)$ we write
\[
	\Aug\dot{\hat{\stateDof}} = \Aug\modelOpDof(\hat{\stateDof}) + \gainFilterDof(\observDof - \Aug\obsOpDof(\hat{\stateDof},t)).
\]
The most classical gain is given by the Extended Kalman Filter (EKF): $\gainFilterDof = \CovDof (\diff \obsOpDof)^\intercal \obsNoiseNormDof$, where $\CovDof$ is obtained from the solution of a Riccati equation on the augmented form
\begin{equation}
	\label{eq:Riccati}
 \dot{\CovDof} = (\diff \Aug\modelOpDof) \CovDof + \CovDof\, (\diff \Aug\modelOpDof)^\intercal - \CovDof (\diff \Aug\obsOpDof)^\intercal \obsNoiseNormDof (\diff \Aug\obsOpDof) \CovDof.
\end{equation} 
The operator $\CovDof$ solution of the Riccati equation is called \emph{covariance} since it can be linked in a stochastic framework to the covariance of the state estimation error evolving during the sequential estimation \cite{Simon06}. When decomposed on the state and parameter
\[
	\CovDof = \begin{pmatrix}
		\stateCovDof & \crossCovDof \\
		(\crossCovDof)^\intercal & \paramCovDof
	\end{pmatrix}
\]
this gives
\[
	\stateGainFilterDof =  \stateCovDof (\diffStateDof  \obsOpDof)^\intercal \obsNoiseNormDof, \quad \paramGainFilterDof =  (\crossCovDof)^\intercal (\diffStateDof  \obsOpDof)^\intercal \obsNoiseNormDof,
\]
and 
\begin{equation}\label{eq:EKF}
	\begin{cases}
				\dotParamCovDof = -(\crossCovDof)^\intercal(\diffStateDof  \obsOpDof)^\intercal \obsNoiseNormDof (\diffStateDof  \obsOpDof) \crossCovDof, \quad \paramCovDof(0) = \paramNoiseCovDof \\	
				\dotCrossCovDof = (\diffStateDof \modelOpDof) \crossCovDof + (\diffParamDof \modelOpDof) \paramCovDof - \stateCovDof (\diffStateDof  \obsOpDof)^\intercal \obsNoiseNormDof (\diffStateDof  \obsOpDof) \crossCovDof, \quad \crossCovDof(0) = 0 \\
				\dotStateCovDof = (\diffStateDof \modelOpDof) \stateCovDof + (\diffParamDof \modelOpDof) (\crossCovDof)^\intercal + \stateCovDof (\diffStateDof \modelOpDof)^\intercal + \crossCovDof (\diffParamDof \modelOpDof)^\intercal \\ \hspace{5.5cm}- \stateCovDof (\diffStateDof  \obsOpDof)^\intercal \obsNoiseNormDof (\diffStateDof  \obsOpDof) \crossCovDof, \quad \stateCovDof(0) = \initNoiseCovDof 
	\end{cases}	
\end{equation}

The practical algorithms --~ for instance time-discrete EKF or UKF for Unscented Kalman Filter --~derived from this formulation are presented in~\ref{append:ekf-ukf}. Even if UKF differs from EKF, it is easy to prove that their analysis which relies on the linearization of their respective estimation error is based on the same error system \cite{MoireauReducedorder}. This is due to the fact that when all operators are assumed to be linear, a finite difference operator or the tangent operator are in fact identical. This implies that we will rely on UKF for the numerical simulation whereas, for the estimator analysis, we will keep EKF which can easily be written in both time-continuous and time-discrete formulations. 

Various sequential estimators can be defined on any dynamical system, and parameter identification can be easily included. This is due to the underlying optimal principles behind the estimator presented before, even if in a non-linear framework we have only relied on approximation rules to extend the optimal filters found in the linear case. However, these estimators have a fundamental drawback since the covariance operator $\CovDof$ is a full matrix of $\M_{\Nstate+\Nparam}(\R)$. Recalling that $\Nstate$ is the dimension of the state variables which discretize the PDEs field variables, the  optimal estimators are intractable for classical finite element models.

\subsection{Reduced-Order Optimal Filters} 
\label{sub:reduced_order_optimal_filters}

A classical first strategy for circumventing the curse of dimensionality consists in assuming a specific reduced-order form for the covariance operators. For example, making the ansatz
\begin{equation}\label{eq:red-order-princ}
	\forall t,\quad \CovDof(t) = \projCovDof(t) \redCovDof(t)^{-1} \projCovDof(t)^\intercal
\end{equation}
with $U$ an invertible matrix of small size $\dimReduce$ and $\projCovDof$ an extension operator, we can show \cite{PM-DC-PLT-08} that with linear operators the solution of the Riccati equation~\eqref{eq:Riccati} in augmented form reduces to
\begin{equation}\label{eq:reducedDynamics}
	\dot{\projCovDof} = \Aug\modelOpDof \projCovDof \text{ and }
	\dot{\redCovDof} = \projCovDof^\intercal \obsOpDof^\intercal \obsNoiseNormDof \obsOpDof \projCovDof,
\end{equation}
which are now computable in practice. In a non-linear framework, the covariance dynamics can then be approximated as in \cite{PM-DC-PLT-08} by extending \eqref{eq:reducedDynamics} as
\begin{equation}\label{eq:nl-reducedDynamics}
	\dot{\projCovDof} = (\diff\, \Aug\modelOpDof) \projCovDof \text{ and }
	\dot{\redCovDof} = \projCovDof^\intercal (\diff \obsOpDof)^\intercal \obsNoiseNormDof (\diff \obsOpDof) \projCovDof.	
\end{equation}
This observer is called \emph{Reduced-Order Extended Kalman Filter} (RoEKF). The time discrete version of this approach is presented in \ref{append:roekf-roukf}, as well as its UKF counterpart. 

\subsection{Luenberger Filters} 
\label{sub:luenberger_filters}

Another way to circumvent the curse of dimensionality is to built a filter which is not based on an underlying optimal criterion. This idea was initially introduced in \cite{Luenberger1971} and therefore is often called \emph{Luenberger} filter and observer. It was quickly popularized in data assimilation \cite{hoke1976initialization,stauffer1990use,Auroux2005} where the curse of dimensionality was very limiting for systems coming from the discretization of PDEs. In data assimilation this strategy is also referred to as the \emph{nudging} approach because the filter is designed to ``gently'' correct the dynamics. 
The goal of the correction is to make the estimation error $\Aug\tilde{\state} = \Aug\state^\target - \Aug\hat{\state}$ converge to 0. In a linear framework, the state error dynamics is 
\[
	\Aug\dot{\tilde{\state}} = (\modelOp - \Aug\gainFilter \obsOp) \Aug\tilde{\state},
\]
thus the strategy is to design a filter gain $\gainFilter$ based on the underlying physics of the problem which makes the dynamics operator $(\modelOp - \gainFilter\obsOp)$ dissipative --~strongly if possible. Such Luenberger filters are designed for the state estimation of specific physical systems and for specific observation operators. For instance, filters have been proposed for transport equations \cite{Auroux2005}, Schr\"{o}dinger problems \cite{bonnabel:hal-00447800},  for waves \cite{Chapelle2011}, beams, elasticity \cite{PM-DC-PLT-09,Xu:2006us} or fluid-structure problems \cite{bertoglio2011}. 
In Section~\ref{sub:luenberger_observer}, this approach will be used for the mechanical problem.

\subsection{Joint state and parameter estimation and coupled system estimation with Luenberger filters} 
\label{sub:joint_state_and_parameter_estimation_or_coupled_system_estimation}

System~\eqref{eq:modelAugNotation} has been defined from the augmented state vector gathering the original state and the parameters. After decomposition on $\state$ and $\param$, this implies that the initial parameters now evolve in time within the observer
\begin{equation}\label{eq:seqEstimator}
	\begin{cases}
		\dot{\hat{\state}} = \modelOp(\hat{\state},\hat{\param},t) + \gainFilter^\state(\observ - \obsOp(\hat{\state})), \quad  \hat{\state} = \stateInit,\\
		\dot{\hat{\param}} = \gainFilter^\param(\observ - \obsOp(\hat{\state})), \quad \hat{\param}(0) = \paramInit,
	\end{cases}
\end{equation}
We have seen that the optimal filter does not make any difference between state and parameters.
However, when dealing with Luenberger observer, there is no direct way to incorporate the parameter identification --~i.e. to design $\gainFilter^\param$ --~from an already defined Luengerger filter $\gainFilter^\state$ on the state. This problem is known as \emph{adaptative filtering} \cite{Baumeister:1987ux,Zhang:2002p137}. In particular in \cite{MoireauReducedorder,PM-DC-PLT-08} a systematic strategy has been proposed to extend a possible Luengerger observer to a joint state and parameters estimation by combining the Luenberger filter and an optimal filter reduced to the remaining parameter space. Taking into account that the parameter space is of much smaller dimension, this strategy allows us to compute a physically-based gain on the large dimensional state and an optimal-based filter on the small parameter space. The resulting observer is proved to converge under adequate assumptions in linear cases and can be extended to non-linear cases. 

When considering multi-physics coupled systems, as the electromechanical problem \eqref{eq:electromech}, the same type of problem occurs as for joint state and parameters estimation. If optimal filters are considered, there is no difficulty to combine them. When measurements are available, and when Luenberger filters have also been designed for the different types of physics, it is also easy to combined them. But, as said before, optimal filters are extremely expensive, and it is not always possible to design a Luenberger filter for all the physical compartments of the problem. It will be justified in Section~\ref{sub:estimation_of_coupled_systems} that, as done for the state and parameter estimation, it is possible for one-way coupled problems to combine a Luenberger filter in one part of the system and a reduced optimal filter on the other part.

\section{Data assimilation for the electromechanical model}
\label{sec:data_assimilation_electromech}

\subsection{State estimation in electrophysiology: a reduced-order approach} 
\label{sub:state_estimation_in_electrophysiology}

The above data assimilation principles can been applied to the system modeling the cardiac electrophysiology. Using a few measures of the electrical potential, the goal is to reduce the uncertainties on the state and the parameters involved in the modeling of cardiac cells. For this system, contrary to the mechanical system that will be presented in the next section, there is no straightforward Luenberger observer. The state has therefore to be filtered differently. 

A possible strategy to overcome the curse of dimensionality induced by the optimal filtering is to discretize the problem on a low-dimensional basis. Here, we propose to build this basis by Proper Orthogonal Decomposition (POD). The POD basis is obtained by keeping the most relevant modes resulting from a Principal Component Analysis of a set of pre-computed solutions. For more details about POD, we refer for example to \cite{kunisch-volkwein-01,rathinam-petzold-04} and to \cite{boulakia-schenone-gerbeau-12,chapelle:hal-00834397} for applications to electrophysiology. 

We denote by $\basisPod \in \M_{\Ndof^\elec,\dimPod}(\R)$ the matrix made of the first $\dimPod$ modes. These modes are orthonormal with respect to a given scalar product, typically in $l^2(\R^{\dimPod})$ or $L^2(\domainHeart)$. Denoting by 
$\massPOD$ the Gramian matrix associated with this scalar product, the POD expansion coefficients of a vector $\stateDof^{\elec}$ are given by $\alpha = \basisPod^\intercal \massPOD \stateDof^{\elec}$. Defining $\stateDof^{\elec}_{\reduced} = \basisPod \alpha$ and $\stateDof^{\elec}_{\perp}=\stateDof^{\elec} - \basisPod \alpha$, the state vector can be decomposed as $\stateDof^{\elec} = \stateDof^{\elec}_{\reduced} +  \stateDof^{\elec}_{\perp}$.

As for the reduced filtering, a first strategy would consist in replacing the full order model by the reduced one: 
\[
	\begin{cases}
		\dot{\stateDof}^\elec_{\reduced} = \modelOpDof^\elec_{\reduced}(\stateDof^\elec_{\reduced},\paramDof^\elec), \\
		\stateDof^\elec_{\reduced}(0) = {\stateInitDof^\elec}_{\reduced} + \initNoise_{{\scriptscriptstyle \stateDof^{\elec}_{\reduced}}}, \\
		\paramDof^\elec = \paramInitDof^\elec + \initNoise_{\paramDof^\elec},
	\end{cases}
\]
where $\modelOpDof^\elec_{\reduced}$ is computed from the finite element matrices $\scalMassDof,\, \laplaceDof_{\rm i}$ and the vectors  $ \vec{\rm F}^\elec,\, \vec{\rm G} $ are projected on the POD space. Then, after discretization, this (small) system can be estimated by  EKF or  UKF (\ref{append:ekf-ukf}). 
%
In doing so, the estimator is known to be stable, but this approach has a drawback: if the POD basis is not rich enough to capture the relevant phenomena, the reduced order dynamics may be a poor approximation of the full order one. Thus we may have a consistency problem in our way of modeling and approximating the system of interest. That is why we prefer another strategy which consists in solving the full order dynamics, and then applying the filter to the projection of the state vector on the POD basis. In other words, we apply RoUKF algorithm (\ref{append:roekf-roukf}) with a reduced variable made of the parameters and the components of the augmented state vector on the POD basis:
\begin{equation}
	\label{eq:split-pod-param}
		\Aug\stateDof^\elec_\dimReduce = \begin{pmatrix}
		\alpha \\
		\paramDof
		\end{pmatrix} \in \R^{\dimReduce} = \R^{\dimPod+\Nparam}
	\text{ such that }
		\Aug\stateDof^\elec = \begin{pmatrix}
			\stateDof^{\elec}_{\perp}\\
			\stateDof^{\elec}_\dimReduce
		\end{pmatrix}	
\end{equation}
Contrary to the first approach, the full model is preserved. The drawback is that it cannot be proved that the part of the error which is not filtered will not grow and then pollute the approximation. Thus, while the first approach could suffer from a consistency problem, the second one can suffer from a stability problem. Nevertheless, in all our test cases, this second strategy proved to be robust and accurate. 

 By combining this strategy with UKF, the following algorithm is obtained to filter the electrical state. The initial condition projector is given by
\[
	\projCovDof(0) = \begin{pmatrix}
		\projCovDof^\perp(0)\\
		\projCovDof^\dimReduce(0)
	\end{pmatrix}
	= 
	\begin{pmatrix}
		\projCovDof^\perp(0) \\
		\begin{bmatrix}
			\projAlphaCovDof(0) \\
			\projCovDof^\paramDof(0)
		\end{bmatrix}
	\end{pmatrix}
	=
	\begin{pmatrix}
		0 \\
		\begin{bmatrix}
			\1 \\
			\1
		\end{bmatrix}
	\end{pmatrix}.
\]
Then we consider an adequate UKF sampling rule composed of weights and particles (see \cite{julier1997nek} or \ref{append:ekf-ukf}), we store the associated weights $(\alpha_i)$ in the diagonal matrix $\weightUKFmat_\alpha$ and precompute specific particles, here the so-called unitary sigma-points (i.e.~with zero mean and unit covariance). More specifically, we consider the $\Nsigma = \dimPod+\Nparam+1$ unitary simplex sigma points, where $\dimPod+\Nparam$ is the reduced space dimension \cite{Julier2002Reduced,MoireauReducedorder}. We denote them by $(I^{i})_{1\leq i \leq \Nsigma}$ and perform at each time step: a first sampling step which generates the particles identified by a subscript $^{[i]}$; a prediction step denoted by an additional superscript $^-$ and a correction step denoted instead by an additional superscript $^+$.
\begin{subequations}\label{eq:ROUKFState}	
	\begin{enumerate}
\item Sampling:
\begin{align}\label{eq:algo_sampling}
		\begin{cases}
			\redCovDofSqrt_n &= \sqrt{(\redCovDof_n)^{-1}} \\[0.1cm]
			\hat{\stateDof}^{[i]+}_{\perp n} &= \hat{\stateDof}_{\perp n}^+ + \projStateCovDofPerp_n \cdot \redCovDofSqrt_n \cdot I^{i},\quad 1\leq i \leq \Nsigma \\[0.1cm]
			\hat{\alpha}^{[i]+}_n &= \hat{\alpha}_n^+ + \projAlphaCovDof_n \cdot \redCovDofSqrt_n^\intercal \cdot I^{i},\quad 1\leq i \leq \Nsigma 	\\[0.1cm]	
			\hat{\paramDof}^{[i]+}_n &= \hat{\paramDof}_n^+ + \projParamCovDof_n \cdot \redCovDofSqrt_n^\intercal \cdot I^{i},\quad 1\leq i \leq \Nsigma 	
		\end{cases}
	\end{align}
\item Prediction:
   \begin{align}\label{eq:algo_prediction}
   	\begin{cases}			
   		\hat{\stateDof}^{[i]-}_{n+1} &= \modelOpDof_{n+1|n}(\hat{\stateDof}^{[i]+}_{n}),\quad 1\leq i \leq \Nsigma \\[0.1cm] 
		\hat{\alpha}^{[i]-}_{n+1} &= \basisPod^\intercal \hat{\stateDof}^{[i]-}_{n+1},\quad 1\leq i \leq \Nsigma \\[0.1cm] %
		\hat{\stateDof}^{[i]-}_{\perp n+1} &=\hat{\stateDof}^{[i]-}_{n+1} - \basisPod \hat{\alpha}^{[i]-}_{n+1},\quad 1\leq i \leq \Nsigma \\[0.1cm] %
   	   	\hat{\stateDof}^-_{\perp n+1} &= \sum_{i=1}^{\Nsigma} \alpha_i \hat{\stateDof}_{\perp n+1}^{[i]-}\\[0.1cm]
   	   	\hat{\alpha}^-_{n+1} &= \sum_{i=1}^{\Nsigma} \alpha_i \hat{\alpha}_{n+1}^{[i]-}\\[0.1cm]
   	   	\hat{\paramDof}^-_{n+1} &= \sum_{i=1}^{\Nsigma} \alpha_i \hat{\paramDof}_{n+1}^{[i]-}
   	\end{cases}
   \end{align}
\item Correction:
   \begin{align}\label{eq:algo_correction}
   	\begin{cases}	
   		\projStateCovDofPerp_{n+1} &= [\hat{\stateDof}_{\perp n+1}^{[*]-}]\weightUKFmat_\alpha [I^{[*]}]^\intercal \\[0.1cm]	
   		\projAlphaCovDof_{n+1} &= [\hat{\alpha}^{[*]-}_{n+1}]\weightUKFmat_\alpha [I^{[*]}]^\intercal \\[0.1cm]	   		
   		\projParamCovDof_{n+1} &= [\hat{\paramDof}^{[*]-}_{n+1}]\weightUKFmat_\alpha [I^{[*]}]^\intercal \\[0.1cm]	   		
   		\observDof_{n+1}^{[i]-} &= \obsOpDof_{n+1}(\hat{\stateDof}^{[i]-}_{n+1}) \\[0.1cm]
		\observDof_{n+1}^{-} &= \sum_{i=1}^{\Nsigma} \alpha_i \observDof_{n+1}^{[i]-} \\[0.1cm]
   		\covObsUKF_{n+1} &=  [\observDof^{[*]-}_{n+1}]\weightUKFmat_\alpha [I^{[*]}]^\intercal \\[0.1cm]
   		\redCovDof^{n+1} &=  \1 + \covObsUKF_{n+1}^\intercal  \obsNoiseCovDof_{n+1}^{-1} \covObsUKF_{n+1} \\[0.1cm]
   		\hat{\stateDof}^{+}_{\perp n+1} &= \hat{\stateDof}^{-}_{\perp n+1} + \projStateCovDofPerp_{n+1} \obsNoiseNormDof_{n+1}  \covObsUKF_{n+1}^\intercal \obsNoiseNormDof_{n+1}  (\observDof_{n+1} - \observDof_{n+1}^{-})\\
   		   		\hat{\alpha}_{n+1}^{+} &= \hat{\alpha}_{n+1}^{-} + \projAlphaCovDof_{n+1} \redCovDof^{n+1}  \covObsUKF_{n+1}^\intercal \obsNoiseNormDof_{n+1}  (\observDof_{n+1} - \observDof_{n+1}^{-})\\
   		   		\hat{\paramDof}_{n+1}^{+} &= \hat{\paramDof}_{n+1}^{-} + \projParamCovDof_{n+1} \redCovDof^{n+1}  \covObsUKF_{n+1}^\intercal \obsNoiseNormDof_{n+1}  (\observDof_{n+1} - \observDof_{n+1}^{-})
   	\end{cases}
   \end{align}
\end{enumerate}
\end{subequations}
where $[I^{[*]}]$ is the matrix concatenating the $(I^{i})$ vectors side by side, and similarly for other vectors \cite{MoireauReducedorder}. From the last three corrections defined in \eqref{eq:ROUKFState}, the state correction reads
\[
\hat{\stateDof}_{n+1}^{+} = \hat{\stateDof}_{n+1}^{-} +\left(\projStateCovDofPerp_{n+1} +\basisPod\projAlphaCovDof_{n+1} \right)\redCovDof^{n+1}  \covObsUKF_{n+1}^\intercal \obsNoiseNormDof_{n+1}  (\observDof_{n+1} - \observDof_{n+1}^{-}).
\]

\subsection{State estimation in mechanics: Luenberger observers} 
\label{sub:luenberger_observer}
To filter the mechanical state, a Luenberger approach is used, following \cite{PM-DC-PLT-09}. Consider the first two equations of \eqref{eq:solid} with
\begin{equation}\label{eq:elastodynamics} 
    \begin{cases}
		\partial_t \disp = \vel, & \quad\mbox{in }\domainHeartRef \\
		\volmass \partial_t \vel  - \Div(\firstPiola)  =  0,& \quad\mbox{in }\domainHeartRef \\
	\end{cases}	
\end{equation}
where the observations are given by
\begin{equation}\label{eq:obselasto}
	\observ^\mech = \disp_{|\domainMeasMech}.
\end{equation}
We then define the estimator by
\begin{equation}\label{eq:elastoLuenberger} 
    \begin{cases}
		\partial_t \hat{\disp} = \hat{\vel} + \gamma \ext_{\domainMeasMech}(\observ^\mech - \hat{\disp}_{|\domainMeasMech}), & \quad\mbox{in }\domainHeartRef \\
		\volmass \partial_t \hat{\vel}  - \Div(\hat{\firstPiola})  =  0,& \quad\mbox{in }\domainHeartRef \\
		\text{same boundary conditions as \eqref{eq:solid}}
	\end{cases}	
\end{equation}
where $\gamma$ is a scalar gain and $\ext_{\domainMeasMech}$ is an extension operator typically given for any displacement field $\vect{d}$ by 
\begin{equation}\label{eq:extension}
	\vect{d}^\text{ext} = \ext_{\domainMeasMech}(\vect{d}) : 
	\begin{cases}
		\Div(\elasticityTensor : \strainLinear(\vect{d}^\text{ext})) = 0 & \quad\mbox{in }\domainHeartRef \backslash \domainMeasMech \\
		\vect{d}^\text{ext} = \vect{d} & \quad\mbox{in } \domainMeasMech \\
		\elasticityTensor : \strainLinear(\vect{d}^\text{ext})) \cdot \normal = k_s \vect{d}^\text{ext}, & \quad \text{ on } \Gamma_{n} \\
		(\elasticityTensor : \strainLinear(\vect{d}^\text{ext})) \cdot \normal = 0, & \quad \text{ on }  \partial \domainHeartRef  \backslash  \Gamma_n
	\end{cases}
\end{equation}
with $\elasticityTensor$ the elasticity tensor coming from the linearization of $\firstPiola$ around $0$ or a given trajectory at a given time $t$.  On a linearized system, it can be proved that the  the state error $\tilde{\state} = \bigl(\begin{smallmatrix}
\tilde{\disp} \\ \tilde{\vel}
\end{smallmatrix}\bigr) = \bigl(\begin{smallmatrix}
\disp - \hat{\disp}\\ \vel -\hat{\vel}
\end{smallmatrix}\bigr) $ between the estimator and the target tends to zero \cite{PM-DC-PLT-09}.

A consistent space discretization is given by
\begin{equation}\label{eq:elastodynamicsDiscrete}
	\begin{cases}
		\dot{\hat{\dispDof}} = \hat{\velDof} + \gamma\extDof(\observDof^\mech - \obsOpDof^\mech(\hat{\stateDof})), \\
		\massDof \ddot{\hat{\dispDof}} + \stressDof(\hat{\dispDof},\hat{\velDof}) = \followPresureDof 
	\end{cases}
\end{equation}
which converges to the real trajectory. Note that the proof of convergence after discretization is in general difficult and should be proved here using the internal viscosity to control the spurious high frequencies introduced by the discretization.

Finally in a state-space form, there exists a Luenberger filter $\gainFilterDof^\mech_{\rm L}$ such that 
\begin{equation}\label{eq:mechStateEstim}
	\dot{\hat{\stateDof}}^\mech = \modelOpDof^\mech(\hat{\stateDof}^\mech) + \gainFilterDof^\mech_{\rm L} (\observDof^\mech - \obsOpDof^\mech(\hat{\stateDof}^\mech)) 
\end{equation}
which converges for any initial error to an observed trajectory. This Luenberger filter, defined in the context of passive non-linear mechanics, was also shown to be robust to the introduction of the active part in the cardiac mechanics \cite{chabiniok-moireau-ea-11,PM-DC-PLT-09}. In this context, the Luenberger filter can be considered as a filter on the reduced space of the displacement and velocity field whereas the internal variables associated with the heart contraction are already stable. Nevertheless, note that the filter robustness  has been demonstrated numerically but remains to be proved theoretically.



\subsection{Estimation of one-way coupled systems} 
\label{sub:estimation_of_coupled_systems}

In this section we present an original strategy to aggregate the already defined filters of each submodel --~namely the electrophysiological model and the mechanical model. This strategy is then justified by an estimation error analysis. To simplify the presentation and analysis of our method, we first consider the following \emph{one way} coupling system with perfectly known parameters. Indeed, we have already seen how a joint state and parameter estimation can be added ``in a second stage'' once the state estimation is proved to be effective. Hence in practice the parameter identification will be considered but here, we can focus on the state as the main difficult aspect of the analysis. We thus consider
\begin{equation}\label{eq:oneWayCouplingState}
	\begin{cases}
		\dot{\stateDof}^\elec = \modelOpDof^\elec(\stateDof^\elec), \quad \stateDof^\elec(0) = \stateInitDof^\elec + \initNoise^\elec \\
		\dot{\stateDof}^\mech = \modelOpDof^\mech(\stateDof^\mech,\stateDof^\elec), \quad \stateDof^\mech(0) = \stateInitDof^\mech + \initNoise^\mech
	\end{cases}
\end{equation}
and we propose to combine the Luenberger state filter on the mechanical part and the reduced order filter of Section~\ref{sub:state_estimation_in_electrophysiology} on the electrical part. We also recall that the observation operator is composed of two concatenated operators
$\obsOpDof^\elec(\stateDof^\elec)$ and $\obsOpDof^\mech(\stateDof^\mech)$ as defined in \eqref{eq:obsOpTotal}.

We can consider a RoEKF without loss by generality at the time-continuous level. Concerning the derivatives, we denote $\diff_\mech$ and $\diff_\elec$ respectively the partial differential with respect to $\stateDof^\mech$ and $\stateDof^\elec$ respectively, whereas we keep $\diff$ when there is no ambiguity or when the differential is total --~\emph{i.e.} with respect to all variables. From the general reduced order formulation of the EKF presented in Section~\ref{sub:reduced_order_optimal_filters} we introduce this time
\[
	\projCovDof = \begin{pmatrix}
		\projCovDof^\elec \\
		\projCovDof^\mech
	\end{pmatrix}
\]
such that
\begin{equation}\label{eq:oneWayCouplingStateFilter}
	\begin{cases}
		\dot{\hat{\stateDof}}^\elec &= \modelOpDof^\elec(\hat{\stateDof}^\elec) + \projCovDof^\elec  \redCovDof^{-1} {\projCovDof}^\intercal {\diff \obsOpDof}^\intercal \obsNoiseNormDof(\observDof - \obsOpDof(\hat{\stateDof})) , \quad \hat{\stateDof}^\elec(0) = \stateInitDof^\elec \\
		\dot{\hat{\stateDof}}^\mech &= \modelOpDof^\mech(\hat{\stateDof}^\mech,\hat{\stateDof}^\elec) + \gainFilterDof^\mech_{\rm L}(\observDof^\mech - \obsOpDof^\mech(\hat{\stateDof}^\mech)) \\
		&\hspace{2cm} +\, \projCovDof^\mech \redCovDof^{-1} {\projCovDof}^\intercal {\diff \obsOpDof}^\intercal \obsNoiseNormDof (\observDof - \obsOpDof(\hat{\stateDof})), \quad \hat{\stateDof}^\mech(0) = \stateInitDof^\mech \\
		\dot{\projCovDof}^\elec &= (\diff \modelOp^\elec(\hat{\stateDof}^\elec)) \projCovDof^\elec, \quad \projCovDof^\elec(0) = \1\\
		\dot{\projCovDof}^\mech &= (\diff_\mech \modelOpDof^\mech(\hat{\stateDof}^\mech,\hat{\stateDof}^\mech) - \gainFilterDof^\mech_{\rm L}\diff\obsOpDof^\mech(\hat{\stateDof}^\mech))\projCovDof^\mech \\
		&\hspace{2cm}+ (\diff_\elec \modelOpDof^\mech(\hat{\stateDof}^\mech,\hat{\stateDof}^\elec))\projCovDof^\elec \quad \projCovDof^\mech(0) = 0\\
		\dot{\redCovDof} &= \projCovDof^\intercal (\diff \obsOpDof)^\intercal \obsNoiseNormDof (\diff \obsOpDof) \projCovDof,\quad \redCovDof(0) = \redCovDof^\elec_\diamond.
	\end{cases}
\end{equation}
Note that in the dynamics \eqref{eq:oneWayCouplingStateFilter} the Luenberger filter only applies on the mechanical part using only the mechanical data. However, the optimal filter strategy allows to benefit from both the electrical and mechanical data to correct the electrophysiology dynamics. This correction is then reverberated to the mechanical dynamics which also reads
\[
	\dot{\hat{\stateDof}}^\mech = \modelOpDof^\mech(\hat{\stateDof}^\mech,\hat{\stateDof}^\elec) + \gainFilterDof^\mech_{\rm L}(\observDof^\mech - \obsOpDof^\mech(\hat{\stateDof}^\mech)) + \projCovDof^\mech (\projCovDof^\elec)^{-1} (\dot{\hat{\stateDof}}^\elec - \modelOpDof^\elec(\hat{\stateDof}^\elec)),
\]
when assuming $\projCovDof^\elec$ invertible. Note that in the particular case where $\modelOpDof^\elec$ is linear the last expression simplifies into
\[
		\dot{\hat{\stateDof}}^\mech = \modelOpDof^\mech(\hat{\stateDof}^\mech,\hat{\stateDof}^\elec) + \gainFilterDof^\mech_{\rm L}(\observDof^\mech - \obsOpDof^\mech(\hat{\stateDof}^\mech)) + \projCovDof^\mech \overbracket[0.5pt][1pt]{(\projCovDof^\elec)^{-1} \hat{\stateDof}^\elec}^{\bm \cdot} ,
\]


\subsubsection{Convergence analysis} 
\label{ssub:convergence_analysis}

To ensure the convergence of the observer, $\projCovDof^\elec$ is assumed to be invertible. Note that $\projCovDof^\elec$ follows the tangent dynamics of the electrophysiological around the estimated trajectory of the model, starting from the initial condition $\projCovDof^\elec(0) = \1$. Therefore, even for a dissipative system, we can assume $\projCovDof^\elec$ to be invertible but potentially ill-conditioned. Moreover the invertibility of $\projCovDof^\elec$ can also be verified numerically in the specific case of concern. We recall that we denote by $\tilde{\stateDof} = \stateDof^\target - \hat{\stateDof}$ the estimation error. The variable $\tilde{\stateDof}$ follows a non-linear dynamics. Then by linearization of this dynamics around the target trajectory we get a linear dynamics satisfied by linearized estimation error denoted by $\errorStateLinDof = 	\bigl( \begin{smallmatrix}
\errorStateLinDof^\elec \\ \errorStateLinDof^\mech 
\end{smallmatrix}\bigr)$. In fact we have
\begin{equation}\label{eq:errorRoEKFcouple}
	\begin{cases}
		\dot{\errorStateLinDof}^\elec &= -\redCovDof^{-1}\projCovDof^\intercal (\diff \obsOpDof)^\intercal \obsNoiseNormDof (\diff \obsOpDof) \errorStateLinDof, \quad \errorStateLinDof^\elec(0) = \initNoise^\elec \\
		\dot{\errorStateLinDof}^\mech &= (\diff_\mech \modelOpDof - \gainFilterDof^\mech_{\rm L}\diff_\mech \obsOpDof^\mech)\errorStateLinDof^\mech + (\diff_\elec \modelOpDof)\errorStateLinDof^\elec + \\&\hspace{3cm} \projCovDof^\mech (\projCovDof^\elec)^{-1} \dot{\errorStateLinDof}^\elec, \quad \errorStateLinDof^\mech(0) = \initNoise^\mech.
	\end{cases}
\end{equation}
Here we introduce the change of variables 
\[
	(\errorStateLinDof^\mech,\errorStateLinDof^\elec) \mapsto (\delta\eta,\delta \mu) = (\errorStateLinDof^\mech - \projCovDof^\mech (\projCovDof^\elec)^{-1} \errorStateLinDof^\elec,(\projCovDof^\elec)^{-1} \errorStateLinDof^\elec),
\]
and obtain
\begin{equation}\label{eq:errorRoEKFchangeVarCouple}
	\begin{cases}
		\dot{\delta\eta} &= (\diff_\mech \modelOpDof - \gainFilterDof^\mech_{\rm L}\diff_\mech \obsOpDof^\mech)\delta\eta, \quad \delta \eta(0) = \initNoise^\mech \\
		\dot{\delta\mu} &= -\redCovDof^{-1}{\projCovDof}^\intercal (\diff \obsOpDof)^\intercal \obsNoiseNormDof (\diff \obsOpDof) \projCovDof \dot{\delta\mu} \\
		&\hspace{1cm} - \redCovDof^{-1}{\projCovDof}^\intercal (\diff_\mech \obsOpDof^\mech)^\intercal \obsNoiseNormDof^\mech (\diff_\mech \obsOpDof^\mech) \delta \eta, \quad \errorStateLinDof^\elec(0) = \initNoise^\elec.
	\end{cases}
\end{equation}
The first equation corresponds to the dynamics of the linearized error studied for the mechanical system. Hence it converges to 0. Therefore the second term in the second equation tends to 0. The homogeneous part of the second equation can then be proved to converge to 0 if the following observability \eqref{eq:observability} condition is satisfied with our linear observation operator --~namely, $\diff \obsOpDof = \obsOpDof$ in our particular example. Namely we expect for all initial error $\errorStateLinDof^\elec(0)$ that
\begin{equation}\label{eq:observability}
		\exists (C,T), \quad \int_0^T \norm{\obsOpDof(\projCovDof (\projCovDof^\elec)^{-1} \errorStateLinDof^\elec)}_\obsNoiseNormDof^2 \geq C \norm{\errorStateLinDof^\elec(0)}^2_{\redCovDof^\elec_\diamond},
\end{equation}
which can be at least verified numerically. In the last observability condition we see that $\projCovDof (\projCovDof^\elec)^{-1} \errorStateLinDof^\elec$ represents the effect of a variation on $\errorStateLinDof^\elec$ on both the electrophysiology and the mechanics. This effect is then observed through $\obsOpDof$. The observability is thus expected to be improved with respect to the situation where only the  electrophysiology is considered. This will be confirmed numerically. 





\subsubsection{Practical algorithm} 
\label{ssub:practical_algorithm}

\begin{figure}[htbp]
	\begin{center}
		
		\begin{tikzpicture}[auto, node distance=1.5cm,>=latex']
		    \node [block] (masterROUKF) {Master\\ RoUKF\\ (Sampling)};		    
		    \node [block, below of=masterROUKF, node distance=3cm] (master) {$[i]$\\Particle \\ Solver};		    
			\node [block, right of=master,node distance=4.6cm] (Mechanics) {$[\Nsigma]$\\Particle  \\ Solver};
			\node [block,left of=master, node distance=4.6cm] (Electric) {$[1]$\\Particle\\ Solver};

\draw [->,line width=1.2pt] (masterROUKF.south) -- node[name=firstsend] {\!\!$\hat{\stateDof}^{[i]+}_{n},\hat{\paramDof}^{[i]+}_{n}$ } (master.north);			
\draw [->,line width=1.2pt] (masterROUKF.south) -- node[name=firstsend] {\!\!\!\!\!\!\! $\hat{\stateDof}^{[\Nsigma]+}_{n},\hat{\paramDof}^{[\Nsigma]+}_{n}$ } (Mechanics.north);			
\draw [->,line width=1.2pt] (masterROUKF.south) -- node[name=firstsend,above] {$\hat{\stateDof}^{[1]+}_{n},\hat{\paramDof}^{[1]+}_{n}$\,\, } (Electric.90);			
 \node [block, below of=master, node distance=3cm] (masterROUKF2) {MASTER \\ RoUKF\\ (Correction)};	
\draw [->,line width=1.2pt] (master.south) -- node[name=firstsend] {\!\!$\hat{\stateDof}^{[i]-}_{n+1},\hat{\paramDof}^{[i]-}_{n+1}$ } (masterROUKF2.north);			
\draw [->,line width=1.2pt] (Mechanics.south) -- node[name=firstsend] {\!\!\!\!\!\!\!\!\!\!\!\!$\hat{\stateDof}^{[\Nsigma]-}_{n+1},\hat{\paramDof}^{[\Nsigma]-}_{n+1}$ } (masterROUKF2.north);			
\draw [->,line width=1.2pt] (Electric.south) -- node[name=firstsend,below] {$\hat{\stateDof}^{[1]-}_{n+1},\hat{\paramDof}^{[1]-}_{n+1}$\,\,\,\,\,\,\,\,\,\, } (masterROUKF2.north);			
\end{tikzpicture}
	
		\caption{Estimation algorithm: ``particle solver'' denotes the electromechanical solver presented in Section~\ref{sec:application_to_a_coupled_electro_mechanical_problem}} 
		\label{fig:masterworkersEstim}
	\end{center}
\end{figure}
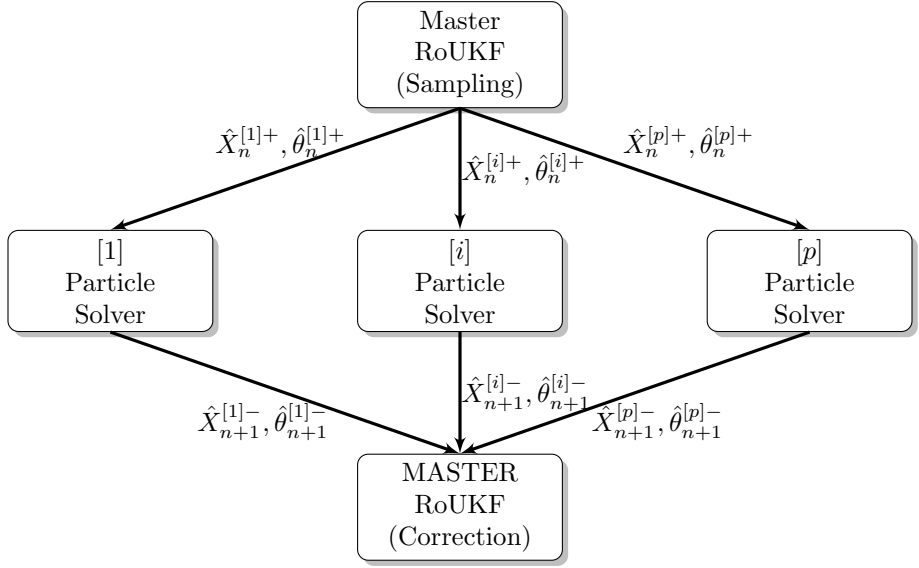

The correction steps of algorithms~\eqref{eq:ROUKF}  (for mechanics) and~\eqref{eq:ROUKFState} (for electrophysiology) are independent of the problem, and can thus be implemented by adopting again a master-slave strategy. In our electromechanical simulator, this yields a two-level interweaved master-slave strategy, the master handling the electromechanical problem becoming a slave of the sequential estimator algorithm. The tasks of the electromechanical master are the data exchange and the correction of unknowns and parameters. The estimation algorithm is summarized in Figure~\ref{fig:masterworkersEstim} for a generic time step. The first data transmission appearing in Figure~\ref{fig:masterworkersEstim} is done in order to update the initial value for each solver. This is handled by the electromechanical master.



\section{Numerical illustrations}
\label{sec:numerical_test_cases}

In this section, the methodology presented in Section~\ref{sub:state_estimation_in_electrophysiology} is applied to the electrophysiological problem of Section~\ref{sub:electrophysiology}. The purpose is twofold: first, assess the reduced order state filtering based on POD; second, illustrate the interest of estimating simultaneaously the state and the parameters.

All the test cases of this Section are performed on the geometry of Figure~\ref{fig:meshes}, including the fibers in the conductivity tensors, and with the parameters typically used to generate healthy ECG \cite{boulakia:inria-00400490}. Synthetic data are generated by applying in \eqref{eq:bidomain2} an external current $I_\text{app}$ during $25\, ms$. Two sets of measurements are considered: (1) 12-lead ECG only or (2) 12-lead ECG enriched with the extra-cellular potential at 8 epicardium points. The first set of measurements is of course the easiest to obtain. The rationale for the second set is to foresee what could be the benefit of including the measurements available from cardiac stimulators or implantable defibrillators. For both cases, these measurements are perturbed by an additive Gaussian white noise of standard deviation of $0.25 mV$ (Figure~\ref{fig:obsNoise}).  Note that another option could be to include endocardial extracellular potentials, if a catheter can provide them. We do not anticipate significant differences with respect to the results obtained with epicardial measurements.

\subsection{ECG based state estimation of the electrophysiological model} 
\label{ssub:electrical_state_results}

\begin{figure}[htbp]
	\begin{center}
		\includegraphics{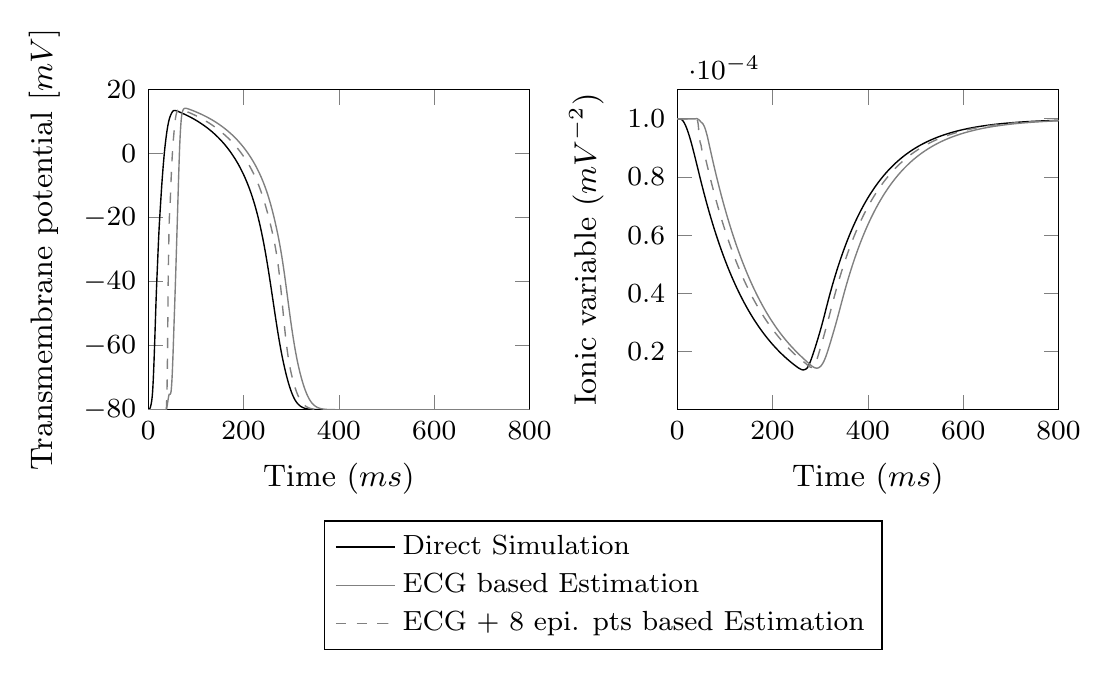}
		\caption{Estimation of the transmembrane potential $\diffV$ (Left) and the ionic variable $\internElec$ (Right). The direct simulation was performed with an external current applied during $25 ms$. The curves represent space-averaged quantities.}\label{fig:Vm_estimation}
		\end{center}
\end{figure}
\begin{figure}[htbp]
\includegraphics[width=\textwidth]{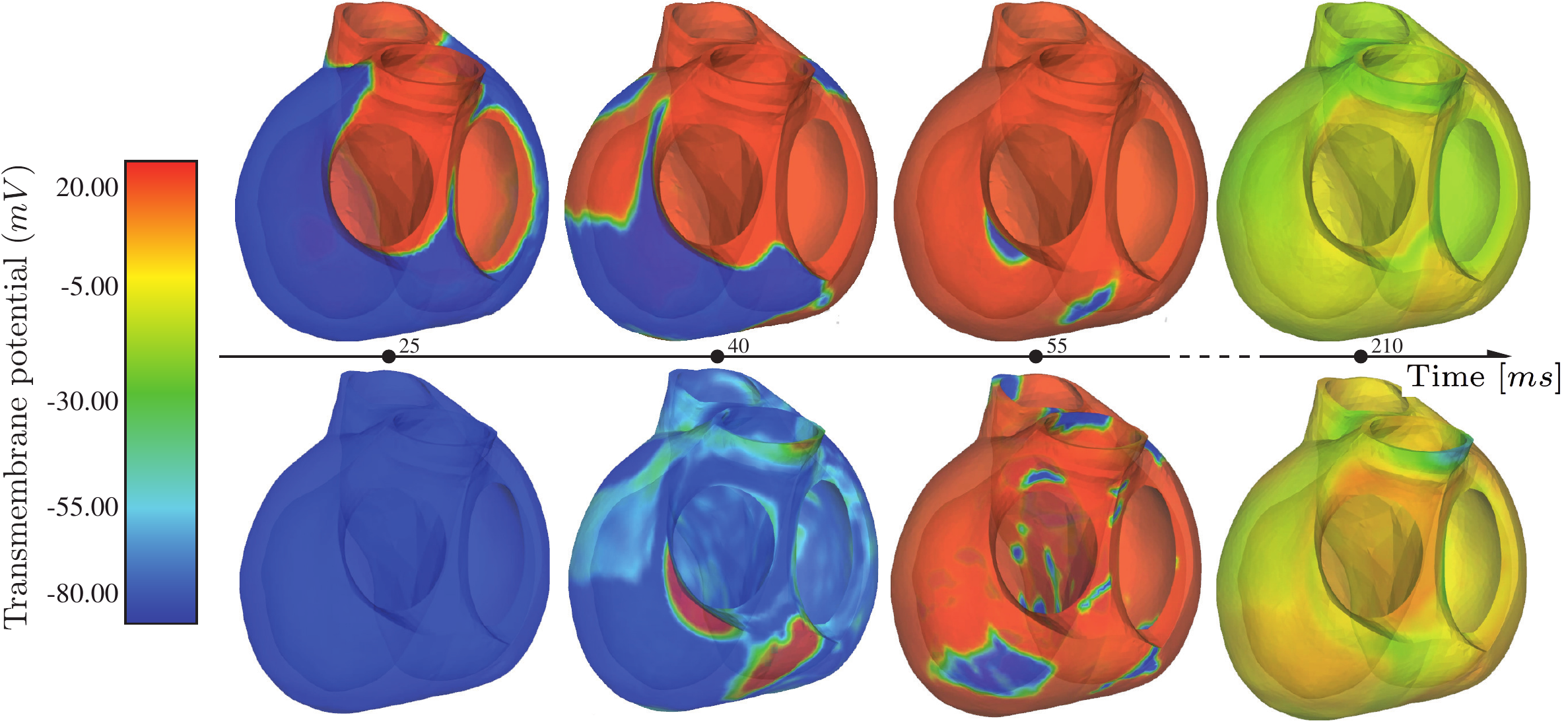}
\caption{Space distribution of the transmembrane potential for the direct simulation (top) and the state estimation with ECG measurements (bottom). There is no initial activation in the bottom simulation in order to generate a strong error in the initial condition. The correction of the filter is applied from $t= 40$ ms. Note the good agreement of both simulation from $t=55 ms$. }\label{fig:vmdistrib}
\end{figure}

Before addressing parameter estimations, a first test is run to assess the reduced order state filter, assuming that all parameters are perfectly known --~for instance with $\tau_{in}=0.8$ and $\tau_{out}=18$. In other words, we study the problem of the state estimation alone applied only tp the electrophysiological model. The POD Basis used in the reduced order state filter \eqref{eq:ROUKFState} was obtained, in this case, from snapshots generated with these parameters, among others. We choose to run the estimation simulation without any external current and compare it with the target simulation starting from $t=40\, ms$. When the estimation starts, the state variables are therefore significantly perturbed with respect to the direct simulation used to generate the synthetic data. In doing so, we introduce a significant state error in our estimation problem. Figure~\ref{fig:Vm_estimation} compares \highlight{the estimated mean value of the transmembrane potential $\diffV$ and the ionic variable $\internElec$} with the target values, while Figure~\ref{fig:vmdistrib} compares the spatial distributions over the whole heart domain. We see that, despite the very large initial error, the estimator succeeds quite well in compensating for the lack of information about the initial stimulation and retrieving in time the accurate electrical state. Indeed, it is remarkable that with a correction only coming from the reduced state filter and applied from $t=40 ms$, the resulting estimated transmembrane potential and ionic current are in good agreement with their target values. Eventually, we point out that, as expected, adding a few measurements on the myocardium slightly improves the result. However this type of measurement clearly requires an invasive procedure \cite{eldar1997transcutaneous}.

\begin{table}[htb]
	\begin{center}
	\begin{tabular}{rcccc}
		\hline Simulation &  $1$    & $2$    & $3$    & $4$ \\ \hline
		$\tau_{in}$       &  $0.8$  & $0.8$  & $1.2$  & $1.2$  \\ \hline
		$\tau_{out}$      &  $18.0$ & $12.0$ & $18.0$ & $12.0$ \\ \hline
	\end{tabular}
	\caption{Set of parameters used to build the POD basis.}
	\label{tab:electrical_results_PODset}
	\end{center}
\end{table}

\begin{figure}[htbp]
	\begin{center}
	\includegraphics{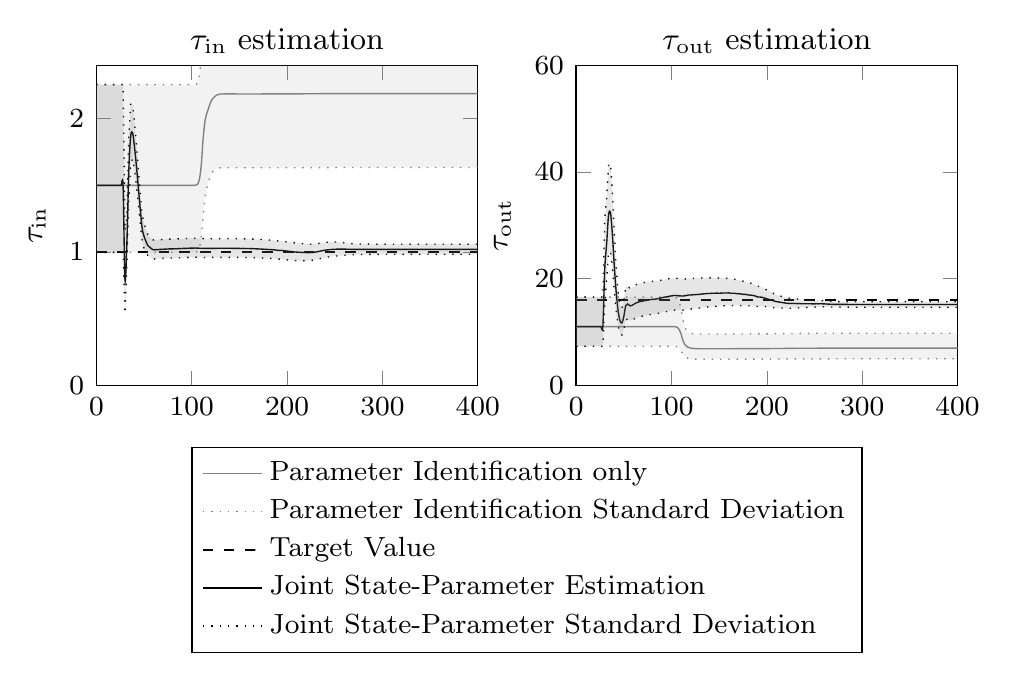}
	\caption{Comparison on $\tin$ (Left) and $\tout$ (Right) estimation, with (black line) and without (gray line) reduced state filtering. These curves clearly show that reducing the uncertainties on the state actually improves the identification of parameters} 	\label{fig:tin1}
	\end{center}	
\end{figure}


\subsection{ECG based parameter identification of the electrophysiological model} 
\label{ssub:electrical_results}
Next, we address the question of the identification of parameters $\tau_{in}$ and $\tau_{out}$ when still only considering the electrophysiological model. The values used to generate the synthetic data are $\tau_{in}=1$ and $\tau_{out}=16$. The initial guess or \emph{a priori} in the inverse problems is $\tau_{in}=1.5,\, \tau_{out}=11$. We point out that neither the values used to generate the synthetic data nor the initial ones are included in the set of solutions used to generate the POD basis. In fact, the POD basis used in the reduced filter~\eqref{eq:ROUKFState} is computed from snapshots obtained with 4 sets of parameters $\tin$ and $\tout$ --~presented in table \ref{tab:electrical_results_PODset} --~which are in a neighbourhood of the expected value, but do not coincide with it as theoretically justified in \cite{chapelle:hal-00834397}. Before launching the estimation, we proceed to a reparametrization of the form $\tin= 2^{\paramDof_1} \tin^\target$ and $\tout=2^{\paramDof_2} \tout^\target$, where $\tin^\target$ and $\tout^\target$ are given and $\paramDof_1$ and $\paramDof_2$ are the new values to be estimated. This is motivated by the fact that the parameters should be maintained positive during the estimation and the uncertainty variance is more naturally centered with respect to a power of 2 of initial parameter. We then compare two strategies: (1) parameter estimation only; (2) joint state-parameter estimation. In strategy (1), the component $\stateDof^\dimReduce$ to be filtered is limited to the parameters $\paramDof$, whereas in strategy (2), $\stateDof^\dimReduce$ includes the parameters $\paramDof$ and the expansion coefficient $\alpha$ on the POD basis (see \eqref{eq:split-pod-param}). 

\begin{table}[h]
	\begin{center}
	\begin{tabular}{ccccc}
		\hline 
		Param. 	  & Target  & \emph{A priori}    & Only Param.     	  & Joint State \&       \\ 
				  & Value   & ($\%$Error)        & Identification     & Parameter Estim. \\ \hline
    	$\tin$    & $1.0$   & $1.5$($50\%$)      & $2.19$($118.8\%$)  & $1.02$($1.99\%$)     \\ \hline
     	$\tout$   & $16.0$  & $11.0$($31.25\%$)  & $6.96$($56.48\%$)  & $15.15$($5.3\%$)     \\ \hline
	\end{tabular}
	\caption{Target vs A priori values for $\tin$ and $\tout$}
	\label{tab:electrical_results_tinout}
	\end{center}
\end{table}

\begin{figure}[htbp]
	\begin{center}
	\includegraphics{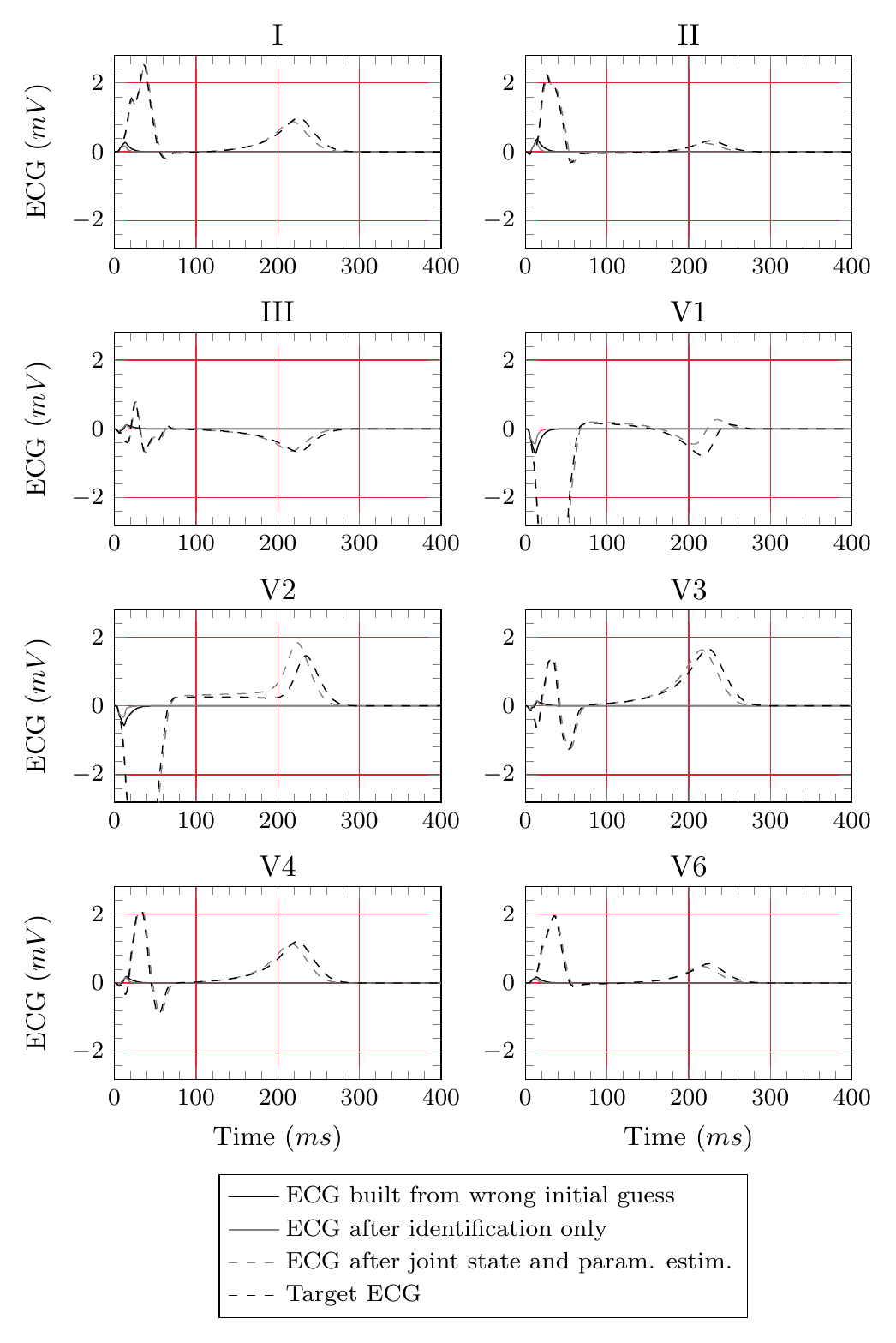}
	\caption{ECGs after state and/or parameter estimation (zoom on a time window of 400 ms)} \label{fig:ecgcomparisonE} 
	\end{center}	
\end{figure}

The evolution of $\tin$ and $\tout$ is reported in Figure~\ref{fig:tin1}, and summarized in Table~\ref{tab:electrical_results_tinout}. In particular we clearly see in Figure~\ref{fig:tin1} the convergence of the estimator on the parameters variable. As an additional illustration of the performance of the joint state-parameter estimator, we plot in Figure~\ref{fig:ecgcomparisonE} 8 leads of the corresponding standard body surface ECGs obtained with the initial guess of the parameters (curve named ``ECG built from wrong initial guess''), with the parameters obtained after a simple parameter estimation (curve named ``ECG after identification only''), with the estimated parameters obtained after a joint state-parameter estimation (curve named ``ECG after joint state and param. estim.''). All these curves should be compared to the one named ``Target ECG'' which corresponds to the measurements used for the estimation. Two comments are in order. First, we note that the ECG obtained with the ``wrong initial guess'' is indeed very different from the target. This means that our initial guess was far from the solution we looked for. Second, we see that the joint state-parameter estimation clearly outperforms the simple parameter estimation (as noted in a different context in \cite{PM-DC-PLT-08}). This means that only a joint state-parameter estimation allows us to produce an electrical state that is really compatible with the observations starting from realistic initial covariance --~namely without forcing \emph{a priori} the estimation to fit the data.

\subsection{Full electromechanical data assimilation} 
\label{sub:results_and_discussion}

In this section, we demonstrate numerically the efficiency of our complete state and parameter estimation chain for coupled models and, more importantly, we show that the estimation of a piecewise constant electrical parameter can be improved by enriching the electrical observations with kinematical observations. The parameter of interest is here $\tclose$, which controls the plateau duration of the action potential of the cardiac cells and is supposed to be heterogeneous in the myocardium. In fact, $\tclose$ is assumed to take different constant values in 4 regions: inner part of the myocardium in the left ventricle, that will be called ``endocardium'' for simplicity; outer part of the myocardium in the left ventricle, that will be called ``epicardium''; a thin region between the endocardium and the epicardium in the left ventricle, that will be called ``M-cell''; and finally the right ventricle. 

Two kinds of synthetic observations are generated from a direct simulation: the electrocardiograms and the displacements of the myocardium. Hence for this numerical illustration, we have at our disposal a $3\D$ field of displacement as it could ideally be processed from a $3\D$ tagged-MRI sequence \cite{SPAMM_MotEst_AKR_SR_2008}. However this is clearly an ideal situation from the mechanical point of view. Nevertheless, we believe that this illustration offers a clear insight into the maximum of information that could be earned with mechanical observations from the point of view of the electrophysiology model.

\begin{figure}[htbp]
	\begin{center}
	\includegraphics{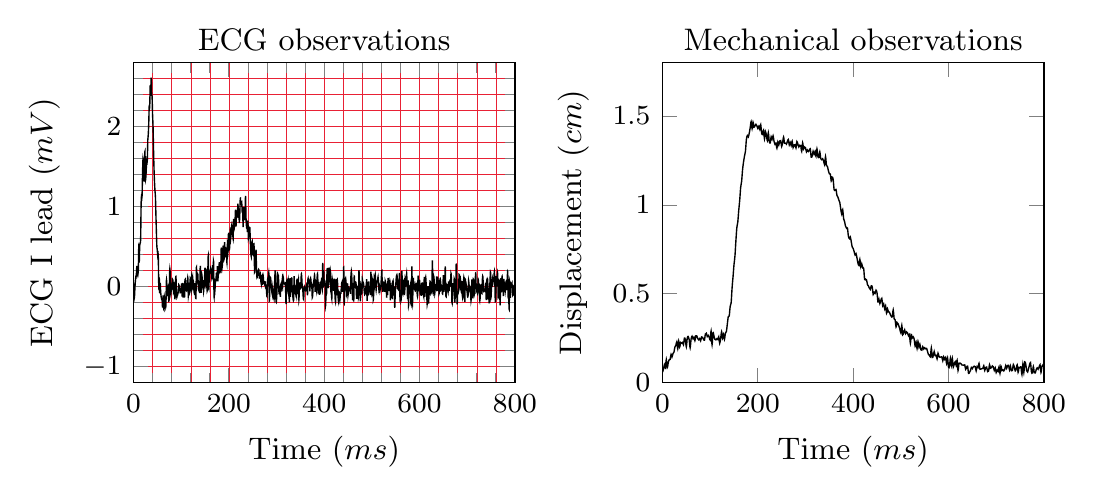}
		\caption{I lead electrocardiogram with additive Gaussian noise (left) and displacement of a mesh point with additive Gaussian noise (right) } \label{fig:obsNoise} 
	\end{center}
\end{figure}

These quantities are perturbed by an additive $10\%$ gaussian white noise in time and space, proportional to the mean amplitude of the signal --~see the corresponding signals in Figure~\ref{fig:obsNoise}. This type of noise is clearly illustrative as it is reasonable to choose independent time observations whereas the spatial distribution of the noise could certainly be more complex than a gaussian noise. In fact, real observations can even be spatially biased. However, this choice of noise is well adapted to evaluate \emph{a minima} the robustness to noise of an estimation method. Facing these noises we follow the recommendation made in \cite{chabiniok-moireau-ea-11} with real data to choose a time-discretized observation norm built from the inverse of the standard deviation $\delta_{\text{obs}}$ by
\[
	\forall n, \obsNoiseNormDof_{n}^\elec = \Delta t_\elec (\delta_{\text{obs}}^\elec)^{-2} \1, \quad \text{ and }, \obsNoiseNormDof_{n}^\mech = \Delta t_\mech (\delta_{\text{obs}}^\mech)^{-2} \1.
\]
with $\delta_{\text{obs}}^\elec \simeq 10\%*2.5\, mV = \, 0.25 mV$ and $\delta_{\text{obs}}^\mech \simeq 10\%*10 \, mm = 1 \,mm$.


%

\begin{table}[t]
	\begin{center}
	\begin{tabular}{rccc}
		\hline Region   & Target value  & Initial guess & Error($\%$) \\ \hline
    	Endocardium     & $140.0$       & $56.0$        & $60\%$      \\ \hline
        M-cell          & $105.0$       & $42.0$        &$60\%$       \\ \hline
       	Epicardium      & $105.0$       & $42.0$        &$60\%$       \\ \hline
   		Right ventricle & $120.0$       & $48.0$        &$60\%$       \\ \hline
	\end{tabular}
	\caption{Target value and initial guess of $\tclose$ in the different regions}
	\label{tab:mech_info_tclose}
	\end{center}
\end{table}

As in Section~\ref{ssub:electrical_results}, the unknowns are reparametrized as $\tclose = 2^\param \tclose^\target$,
where $\tclose^\target$ corresponds to the initial guess. The values to be estimated and the initial guess are gathered in Table~\ref{tab:mech_info_tclose}. 

\begin{figure}[htbp]
	\begin{center}
	\includegraphics{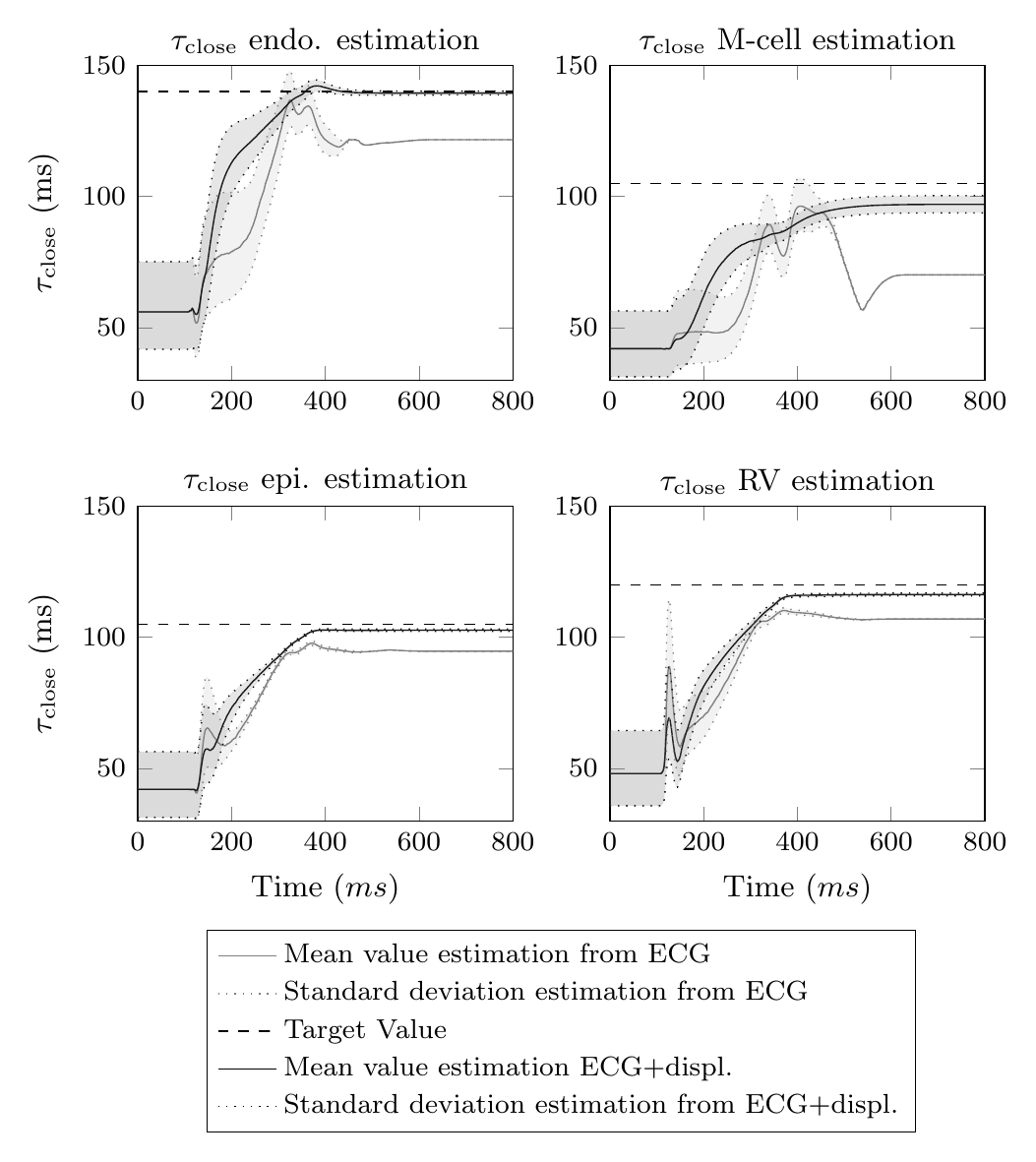}
	\caption{$\tclose$ estimation on endocardium, m-cell, epicardium and right-ventricle regions, comparison ECG observations vs ECG+mechanical observations}\label{fig:tcloseestim}
	\end{center}	
\end{figure}

\begin{figure}[htbp]
	\begin{center}
	\includegraphics{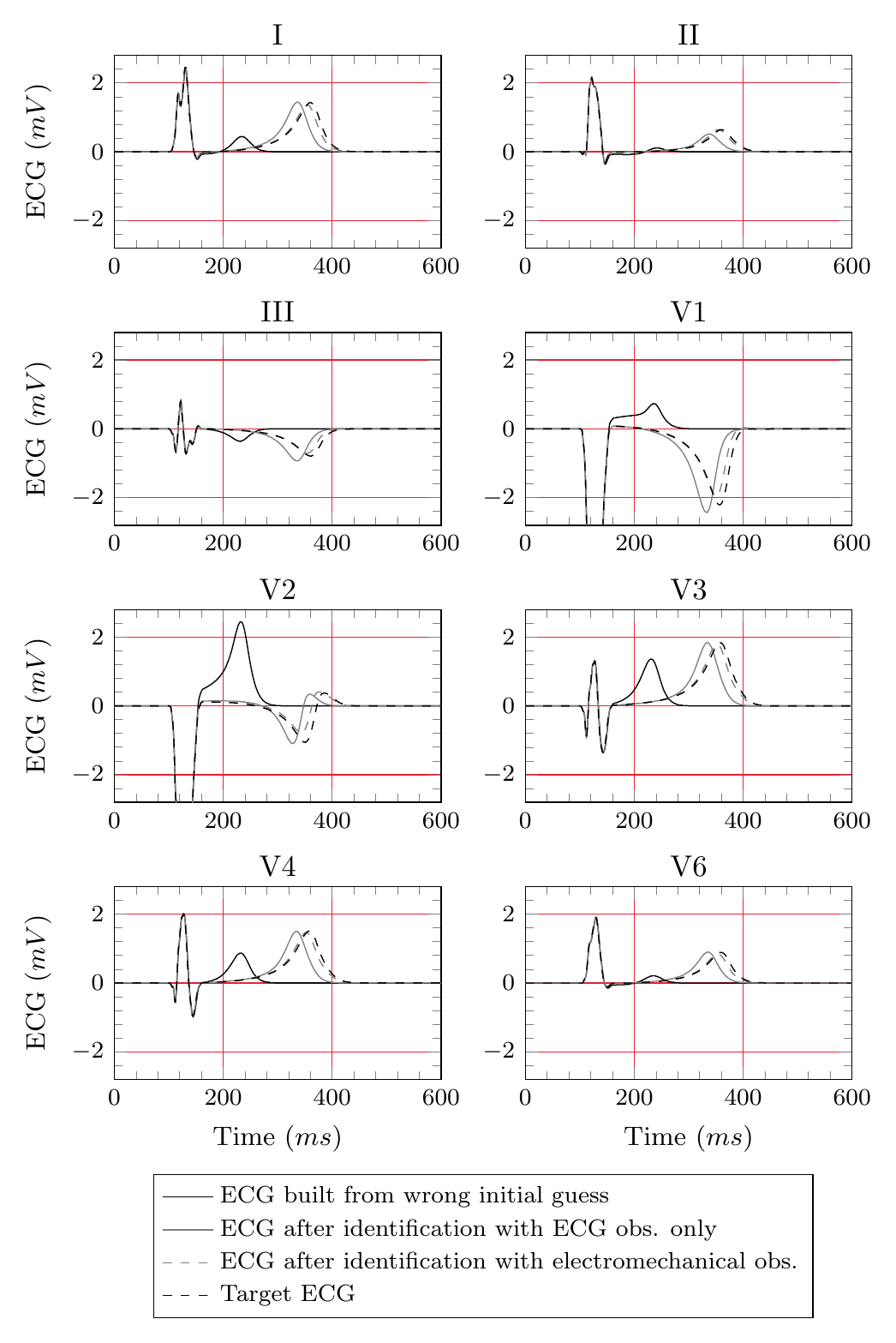}
	\caption{ECGs obtained before and after estimation, compared to the ECG used as a measurement} \label{fig:ecgcomparisonEM} 
	\end{center}	
\end{figure}

Figure~\ref{fig:tcloseestim} shows the time evolution of $\tclose$ during the sequential estimation in each region for electrical measurements only and for the electromechanical --~ECG plus displacements --~available measurements. The benefit of the electromechanical measurements is striking, especially in the M-cell region. Then we present the final estimated values and the error with respect to the target values in Table~\ref{tab:EM_results_estim_tclose_E}. We also report the final estimated values for electromechanical measurements after three heart beats where the estimation error is even more reduced. With the electrical measurements only, we were not able to run three heart beats, because the estimated parameters diverged too much from physiological values. On the contrary, with electromechanical measurements, the results keep improving along the three heart beats. Thus, it clearly appears that the estimation is much more accurate and robust when the electromechanical measurements are taken into account. Our data assimilation procedure results can also be evaluated in the light of ECGs that the model can produce. Indeed, we plot in Figure~\ref{fig:ecgcomparisonEM} the reference ECG (in dashed black), the ECG corresponding to the initial guess (in black), the ECG obtained from the estimated values with the electrical measurement (in gray) and finally the ECG obtained from the estimated values with the electromechanical measurements (in dashed gray). We clearly see that the ECG corresponding to the parameter initial guess has a completely false T-wave, whereas after the estimation procedure the ECGs are much closer to the reference, especially with electromechanical measurements. 

\begin{figure}
	\begin{center}
	\includegraphics{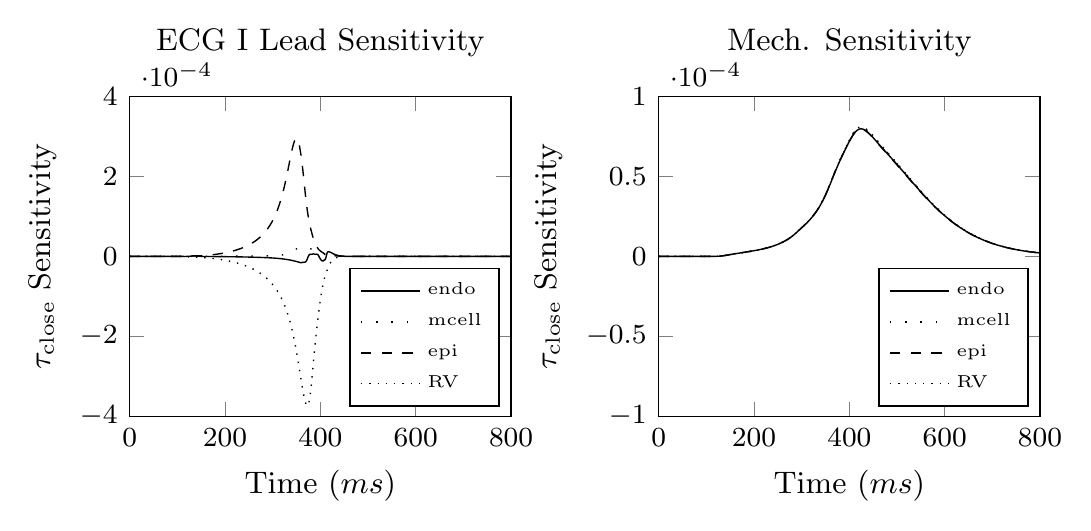}
	\caption{Measurements sensitivity with respect to $\tclose$: ECG Lead I sensitivity (left) and mechanical sensitivity (right, the four curves are superimposed)}
	\label{fig:tclsens}
	\end{center}
\end{figure}
In order to quantitatively understand the relevance of using the mechanical observations to estimate an electrical parameter, a sensitivity analysis is performed. To simplify the presentation of this last calculus, let us consider a time continuous ROEKF for the parameter estimation. We recall that the sensitivity of the model with respect to its parameters is directly estimated with $\projStateCovDof$ in the ROEKF filter \cite{PM-DC-PLT-08}. In our practical case we recall that $\projStateCovDof$ is a matrix which consists of 4 columns vectors (one for each region of the Table~\ref{tab:mech_info_tclose}) of dimension the number of the electrical model degrees of freedom plus the mechanical model degrees of freedom. Therefore, the measurements sensitivity with respect to the parameters is estimated by $\obsOpDof \projStateCovDof$. Note that we have already seen such sensitivity terms in the observability condition \eqref{eq:observability} introduced in the convergence study of our coupled estimator. In our practical case, we focus on the sensitivity $s^\elec$, resp. $s^\mech$, of the first lead of the ECG, resp. the mean value of the measured displacements norm in the myocardium, with respect to $\tclose$, and we consider normalized quantities. We denote by $(\observDof^\elec)_1$ the first index of $\observDof^\elec$ where the first lead of the ECG is gathered. We thus define 
\[
	s^\elec(t) = \frac{{\tclose}_\diamond}{\bar{\observ}^e_1} (\obsOpDof^\elec \projCovDof^\elec(t))_1, \text{ with } \bar{\observ}^\elec_1 = \frac{1}T \!\! \int_0^T \!\!\! \abs{(\observDof^\elec(t))_1} \, dt,
\]
and $T$ is typically the heart beat duration. Identically we define
\[
	s^\mech(t) = \frac{{\tclose}_\diamond}{\bar{\observ}^m} ((\obsOpDof^\mech \projCovDof^\mech(t))^\intercal \obsNoiseNormDof^\mech \obsOpDof^\mech \projCovDof^\mech(t))^{\frac{1}2}, \text{ with } \bar{\observ}^\mech = \frac{1}T \!\! \int_0^T \!\!\! ((\observDof^\mech)^\intercal \obsNoiseNormDof^\mech \observDof^\mech)^{\frac{1}2} \, dt.
\]
In practice when using the ROUKF the computations are very similar but time integrations is replaced by time iteration summations whereas $\obsOpDof \projCovDof(t)$ are replaced by $\covObsUKF_{n+1}(\projParamCovDof_{n+1})^{-1}$ computed from \eqref{eq:ROUKFState}. Our two normalized sensitivities are represented in Figure~\ref{fig:tclsens} for the four different regions. In the epicardium and in the right ventricle, the sensitivity of the electrical measurement is higher than the mechanical one, as expected. But interestingly, the electrical sensitivity is lower than the mechanical one in the M-cell and endocardium regions. This is particularly the case here where we observe the displacement in the whole heart. However recent advances in tagged-MRI allows to expect that intra-myocardial kinematics will be available in the future whereas it is out of reach for non-invasive electrical measurement. In addition, the mechanical observation is affected by $\tclose$ over a longer time window. Hence, the mechanical measurements therefore increase the time interval during which the observer has the opportunity to correct the value of $\tclose$.

\begin{table}
		\begin{tabular}{rccc}
			\hline  Parameter        & ECG only Estim       & ECG+Mech Estim      & ECG+Mech Estim      \\
			                         & after 1 beat         & after 1 beat        & after 3 beat        \\
		                             & (Error $\%$)         & (Error $\%$)		  & (Error $\%$)		\\\hline
			$\tclose^{\text{endo}}$  & 121.63 ($13.12\% $)  & 139.42  ($0.41\% $) & 140.9   ($0.64\% $) \\\hline
			$\tclose^{\text{mcell}}$ & 70.12  ($33.22\% $)  & 96.98   ($7.64\% $) & 104.23  ($0.73\% $) \\\hline
			$\tclose^{\text{epi}}$   & 94.64  ($9.87\% $)   & 102.68  ($2.21\% $) & 104.28  ($0.69\% $) \\\hline
			$\tclose^{\scriptscriptstyle \text{RV}}$ 
									 & 106.94 ($10.88\% $)  & 116.26  ($3.12\% $) & 118.74  ($1.05\% $) \\\hline
		\end{tabular}
		\caption{Identification of $\tclose$ in the various regions and for the different scenarios}
		\label{tab:EM_results_estim_tclose_E}
\end{table}

%

In conclusion, this last example is a rich illustration of our sequential strategy, with a state filter based on a Luenberger approach for the mechanics and on a reduced order optimal filter for the electrophysiology and the parameters. In particular it allows us to numerically justify the importance of having a correction for every type of uncertainties, especially for state uncertainties. Moreover we numerically demonstrate the interest of taking advantage of the multi-physics nature of a problem to identify some parameters --~here $\tclose$ accounting for the plateau duration --~which have an impact on the coupled physics.

\section{Conclusion} 
\label{sec:conclusion_}
In this work, a complete framework for the joint state and parameter estimation of a complex multi-physics problem has been presented. The strategy consists of Luenberger observers and reduced order optimal observers. The robustness of the approach is ensured by the fact that each component of this system is corrected by the data. In particular, the importance of state estimation when performing a parameter identification has been illustrated. For weakly coupled systems we have justified that our filter combination and aggregation lead to a reduction of the estimation error. It is expected that the same type of results could be proved for fully coupled problems in the future. This general strategy has been applied to an electromechanical model of the heart but other coupled systems of interest in various fields of application can be estimated with this approach. In this framework, we have also shown that the use of mechanical data could considerably improve the observability of the problem.

The direct electromechanical model is clearly more demanding than the direct electrophysiology model. However, we have shown that the associated inverse problem which consists in identifying its parameters from all the available data --~mechanical and electrical --~is better defined than the pure electrophysiology inverse problem.

In conclusion, this work provides a promising  new strategy to address the classical inverse problem of electrocardiography, and suggests a new way to reduce its well-known ill-posedness. The numerical results, based on synthetic data, are presented primarily to validate the methodology. This validation being successfully achieved, an important  issue to be addressed in future works concerns the adequacy and accuracy of the models. Another important perspective is the validation against real data. Using real data  will necessitate a complete patient workflow where  anatomical data, ECG and Tagged-MRI (or at least Cine-MRI based on the results of \cite{chabiniok-moireau-ea-11}) are acquired and post-processed. In this respect, we point out that our strategy is indeed quite demanding in term of data acquisition and processing. This must be ultimately compared --~in term of estimation robustness and accuracy versus difficulty--~to existing methods such as electrocardiographic imaging (ECGI) \cite{ramanathan-ghanem-jia-04} where only ECG-like measurements are necessary with special vests recording up to 224 body-surface potentials.


\section*{Acknowledgement} 

The authors would like to thank the anonymous reviewers for their valuable comments and suggestions.

\appendix

\section{Time-discrete EKF and UKF filters}
\label{append:ekf-ukf}

The time-discretization of the optimal estimator presented in Section~\ref{sub:the_extended_kalman_filter} is based on the principle that the optimality should be conserved at the discrete level. Let us first denote a stable and consistent discretization of the original model by
\begin{equation}\label{eq:discreteModel}
	\begin{cases}
		\stateDof_{n+1} = \modelOpDof_{n+1|n}(\stateDof_{n},\theta_{n}) \\
		\stateDof_0 = \stateInitDof + \initNoise^\state \\
		\paramDof = \paramInitDof + \initNoise^\param
	\end{cases}
\end{equation}
and define a discrete-time functional 
\[
	\mathscr{J}_N(\initNoiseStateDof,\initNoiseParamDof) = \frac{1}2 \norm{\initNoiseStateDof}^2_{\initNoiseCov^{-1}} + \frac{1}2 \norm{\initNoiseParamDof}_{\paramNoiseCov^{-1}}^2 +\frac{1}2 \sum_{k=0}^{N_{T}} \norm{\observDof_k - \obsOpDof(\stateDof_{k,[\initNoiseStateDof,\initNoiseParamDof]})}_{\obsNoiseNorm_k}^2 
\]
with $\obsNoiseNorm_k = \Delta t \obsNoiseNorm$ when a fixed time-step is considered. This choice of discrete observation norm $\obsNoiseNorm_k$ ensures that the discrete-time functional $\mathscr{J}_N$ is consistent with respect to the continuous-time functional $\mathscr{J}$. The discrete-time optimal filter is then defined by $\hat{\stateDof}_n = \stateDof_{n,[\argmin \mathscr{J}_n]}$. The discrete-time counterpart of \eqref{eq:EKF} for the optimal linear estimator can be deduced with a prediction-correction time scheme. We denote by a superscript $-$ the prediction computations and by a superscript $+$ the correction computations.
\begin{subequations}\label{eq:KalmanDiscrete}	
	\begin{enumerate}
	\item Prediction: 
   \begin{align} 
   	\begin{cases}			
   		\hat{\stateDof}_{n+1}^- &= \modelOpDof_{n+1|n}(\hat{\stateDof}_n^+,\hat{\paramDof}_n^+) \\
   		\hat{\paramDof}_{n+1}^- &= \hat{\paramDof}_n^+ \\
		\stateCovDofPred_{n+1} &= (\diffStateDof \modelOpDof_{n+1|n}) \stateCovDofCorr_n (\diffStateDof \modelOpDof_{n+1|n})^\intercal  \\ 
		& \hspace{1cm}  +\, (\diffStateDof \modelOpDof_{n+1|n}) \crossCovDofCorr_n (\diffParamDof \modelOpDof_{n+1|n})^\intercal \\
        & \hspace{1cm} 	+\, (\diffParamDof \modelOpDof_{n+1|n}) (\crossCovDofCorr_n)^\intercal (\diffStateDof \modelOpDof_{n+1|n})^\intercal  \\ 
		& \hspace{1cm} 	+\, (\diffParamDof \modelOpDof_{n+1|n}) \paramCovDofCorr_n (\diffParamDof \modelOpDof_{n+1|n})^\intercal \\
		\crossCovDofPred_{n+1} &= (\diffStateDof \modelOpDof_{n+1|n}) \crossCovDofCorr_n   + 
                                  (\diffParamDof \modelOpDof_{n+1|n}) \paramCovDofCorr_n \\
		\paramCovDofPred_{n+1} &= \paramCovDofCorr_{n}
   	\end{cases}
   \end{align} 
	\item Correction:
   \begin{align}\label{eq:KalmanCorrection}
   	\begin{cases}	
   		\CovDof^+_{n+1} &= \Bigl((\diff \obsOpDof_{n+1})^\intercal \obsNoiseNormDof_{n+1} (\diff \obsOpDof_{n+1}) + (\CovDof^-_{n+1})^{-1}\Bigr)^{-1} \\
		\stateGainFilterDof_{n+1} &= \stateCovDofCorr_{n+1} (\diffStateDof \obsOpDof_{n+1})^\intercal \obsNoiseNormDof_{n+1} \\
		\paramGainFilterDof_{n+1} &= (\crossCovDofCorr_{n+1})^\intercal (\diffStateDof \obsOpDof_{n+1})^\intercal \obsNoiseNormDof_{n+1} \\
		\hat{\stateDof}_{n+1}^+  &= \hat{\stateDof}_{n+1}^- + \stateGainFilterDof_{n+1}(\observDof_{n+1}- \obsOpDof_{n+1}(\hat{\stateDof}_{n+1}^-)) \\
		\hat{\paramDof}_{n+1}^+  &= \hat{\paramDof}_{n+1}^- + \paramGainFilterDof_{n+1}(\observDof_{n+1}- \obsOpDof_{n+1}(\hat{\stateDof}_{n+1}^-))
   	\end{cases}
   \end{align}
\end{enumerate}
\end{subequations} 
This algorithm can be interpreted as prediction-correction time discretization of \eqref{eq:EKF}.

Another way to propose an efficient estimator is to replace at the time-discrete level the tangent operators by a finite difference interpolation based on the computation of the original operator on numerous sampling points. This is the case for Ensemble Kalman Filter (EnKF) \cite{evensen2009data} or the Unscented Kalman Filter (UKF) \cite{julier1997nek}. In this article we focus on the second one which is a discrete-time estimator based on sampling particles called \emph{sigma-points} helping to replace the tangent computations. Let us introduce the so-called \emph{unitary sampling points} $I^{[i]}$ and weight $\alpha_i$ with the following rules
\begin{equation}
		\sum_{i=1}^{\Nsigma} \alpha_i I^{[i]} = 0, \quad
		\sum_{i=1}^{\Nsigma} \alpha_i I^{[i]} \cdot  {I^{[i]}}^\intercal = \1,
\end{equation}
so that, at each time step, the sigma-points can be generated around the estimated values based on the covariance estimation by 
\begin{subequations}\label{eq:UKFDiscrete}
\begin{equation}
	\begin{pmatrix}
		\hat{\stateDof}_{n+1}^{[i]+} \\ \hat{\paramDof}_{n+1}^{[i]+}
	\end{pmatrix} = 	
	\begin{pmatrix}
		\hat{\stateDof}_{n+1}^+ \\ \hat{\paramDof}_{n+1}^+
	\end{pmatrix}+ \sqrt{\CovDof^+_{n+1}} I^{[i]},
\end{equation}
then we compute the prediction
\begin{equation}
	\begin{cases}
		\hat{\stateDof}_{n+1}^{[i]-} = \modelOpDof_{n+1|n}(\hat{\stateDof}_{n}^{[i]+},\hat{\paramDof}_{n}^{[i]+}), \quad &\hat{\stateDof}_{n+1}^- = \sum_{i=1}^{\Nsigma} \alpha_i \hat{\stateDof}_{n+1}^{[i]-}, \\
		\hat{\paramDof}_{n+1}^{[i]-} \,\, = \hat{\paramDof}_{n}^{[i]+}, \quad & \hat{\paramDof}_{n+1}^- = \sum_{i=1}^{\Nsigma} \alpha_i \hat{\paramDof}_{n+1}^{[i]-} = \hat{\paramDof}_{n}^+,
	\end{cases}
\end{equation}
and the corresponding observations
\begin{equation}
	\hat{\observDof}_{n+1}^{[i]-} = \obsOpDof(\hat{\stateDof}_{n+1}^{[i]-}),\quad  \hat{\observDof}_{n+1}^{-} = \sum_{i=1}^{\Nsigma} \alpha_i \hat{\observDof}_{n+1}^{[i]}.
\end{equation} 
The gain is then defined by
\begin{equation}
\begin{cases}
		\crossObsCovDof_{n+1} = \sum_{i=1}^{\Nsigma} \alpha_i(\hat{\stateDof}_{n+1}^{[i]-}-\hat{\stateDof}_{n+1}^-)(\hat{\observDof}_{n+1}^{[i]-}-\hat{\observDof}_{n+1}^{-})^\intercal,\\
		\crossObsCovParamDof_{n+1} = \sum_{i=1}^{\Nsigma} \alpha_i(\hat{\paramDof}_{n+1}^{[i]-}-\hat{\paramDof}_{n+1}^-)(\hat{\observDof}_{n+1}^{[i]-}-\hat{\observDof}_{n+1}^{-})^\intercal, \\
		\obsCovDof_{n+1} = \sum_{i=1}^{\Nsigma} \alpha_i(\hat{\observDof}_{n+1}^{[i]-}-\hat{\observDof}_{n+1}^{-})(\hat{\observDof}_{n+1}^{[i]-}-\hat{\observDof}_{n+1}^{-})^\intercal + \obsNoiseNorm_{n+1}^{-1}, \\
		\stateGainFilterDof_{n+1} = (\crossObsCovDof_{n+1}) \cdot (\obsCovDof_{n+1})^{-1},\\
		\paramGainFilterDof_{n+1}=(\crossObsCovParamDof_{n+1}) \cdot (\obsCovDof_{n+1})^{-1}
\end{cases}
\end{equation}
so that, with the covariance
\[
	\Aug\CovDof_{n+1}^- = \sum_{i=1}^{\Nsigma} \alpha_i \begin{pmatrix}
		\hat{\stateDof}_{n+1}^{[i]-} \\ \hat{\paramDof}_{n+1}^{[i]-}
	\end{pmatrix} \begin{pmatrix}
		\hat{\stateDof}_{n+1}^{[i]-} \\ \hat{\paramDof}_{n+1}^{[i]-}
	\end{pmatrix}^\intercal, \quad \Aug\crossObsCovDof_{n+1} = \begin{pmatrix}
		 \crossObsCovDof_{n+1} \\ \crossObsCovParamDof_{n+1}
	\end{pmatrix},
\]
we still have
\begin{equation}
\begin{cases}
	\hat{\stateDof}_{n+1}^+  = \hat{\stateDof}_{n+1}^- + \stateGainFilterDof_{n+1}(\observDof_{n+1}- \hat{\observDof}_{n+1}^{-}), \\
\hat{\paramDof}_{n+1}^+  = \hat{\paramDof}_{n+1}^- + \paramGainFilterDof_{n+1}(\observDof_{n+1} - \hat{\observDof}_{n+1}^{-}), \\	
	\Aug\CovDof_{n+1}^+ = \Aug\CovDof_{n+1}^- - \Aug\crossObsCovDof_{n+1}  (\obsCovDof_{n+1})^{-1} 	(\Aug\crossObsCovDof_{n+1})^\intercal.
\end{cases}
\end{equation}
\end{subequations}

This filter has the great advantage to be much easier to compute numerically since, contrary to EKF, the tangent operators computations are no longer required. 

\newpage
\section{Reduced-Order filters: time discrete RoEKF and RoUKF}
\label{append:roekf-roukf}

The reduced-order filtering concept presented in Section~\ref{sub:reduced_order_optimal_filters} can be directly applied to the time and space discretized versions of the equations, leading to a time-discrete RoEKF:
\begin{subequations}\label{eq:ROKalmanDiscrete}	
	\begin{enumerate}
	\item Prediction:
   \begin{align} 
   	\begin{cases}			
   		\Aug\hat{\stateDof}_{n+1}^- &= \Aug\modelOpDof_{n+1|n} (\Aug\stateDof_{n}^+) \\
   		\projCovDof_{n+1} &= (\diff\, \Aug\modelOpDof_{n+1|n})  \projCovDof_n
   	\end{cases}
   \end{align}
	\item Correction:
   \begin{align}
   	\begin{cases}	
   		\redCovDof_{n+1} &= \redCovDof_{n} + \projCovDof_{n+1}^\intercal (\diff  \obsOpDof_{n+1})^\intercal \obsNoiseNormDof_{n+1} (\diff \obsOpDof_{n+1}) \projCovDof_{n+1} \\
		\gainFilterDof_{n+1} &=  \projCovDof_{n+1} \redCovDof_{n+1}^{-1} \projCovDof_{n+1}^\intercal (\diff  \obsOpDof_{n+1})^\intercal\obsNoiseNormDof_{n+1} \\
		\Aug\hat{\stateDof}_{n+1}^+  &= \Aug\hat{\stateDof}_{n+1}^- + \gainFilterDof_{n+1}(\observDof_{n+1} - \obsOpDof_{n+1}(\Aug\hat{\stateDof}_{n+1}^-))
   	\end{cases}
   \end{align}
\end{enumerate}
\end{subequations}

This  algorithm has been applied in \cite{PM-DC-PLT-08} for parameter identification by reducing the uncertainty space to the parameter space. Indeed, the extension $\projCovDof$ is initially decomposed into
\[
	\projCovDof(0) = \begin{pmatrix}
		\projStateCovDof \\
		\projParamCovDof
	\end{pmatrix}
	= \begin{pmatrix}
		0 \\
		\1
	\end{pmatrix}.
\]
Using the fact that the parameters dynamics is null, we can easily proved that 
\begin{equation}\label{eq:Ltheta}
	\forall t>0,\quad \projParamCovDof = \1.
\end{equation}
Then, by direct computation, we have the continuous-time formulation of the RoEKF in a parameter identification context 
\begin{equation}\label{eq:paramRoEKF}
	\begin{cases}
		\dot{\hat{\stateDof}} = \modelOpDof(\hat{\stateDof},\hat{\paramDof}) + \projStateCovDof \dot{\hat{\paramDof}}, \quad \hat{\stateDof}(0) = \stateInitDof, \\
		\dot{\hat{\paramDof}} = \redCovDof^{-1}{\projStateCovDof}^\intercal (\diffStateDof \obsOpDof)^\intercal \obsNoiseNormDof (\observDof - \obsOpDof(\hat{\stateDof})), \quad \hat{\paramDof}(0) = \paramInitDof, \\
		\dot{\projStateCovDof} = (\diffStateDof \modelOpDof) \projStateCovDof + \diffParamDof \modelOpDof,\quad \projStateCovDof(0) = 0, \\
		\dot{\redCovDof} = {\projStateCovDof}^\intercal (\diffStateDof \obsOpDof)^\intercal \obsNoiseNormDof (\diffStateDof \obsOpDof) \projStateCovDof, \quad \redCovDof(0) = \redCovDof_*.
	\end{cases}
\end{equation}
In~\eqref{eq:paramRoEKF}, we see that the filter correction appears in the parameters dynamics. Then, the correction on the parameters is reverberated to the state by the mean of $\projStateCovDof$ which can be interpreted from its proper dynamics as the sensitivity of the model with respect to the parameters. Therefore even for a reduced order strategy on the parameters the global filter corrects the two components of the model, namely the parameters but also the state.

The same principles can be applied after the time discretization of the model to get
\begin{subequations}\label{eq:paramRoEKFDiscrete}	
	\begin{enumerate}
	\item Prediction:
   \begin{align} 
   	\begin{cases}			
   		\hat{\stateDof}_{n+1}^- &= \modelOpDof_{n+1|n} (\hat{\stateDof}_{n}^+,\hat{\paramDof}_{n}^+) \\
		\hat{\paramDof}_{n+1}^- &= \hat{\paramDof}_{n}^+ \\
   		\projStateCovDof_{n+1} &= (\diffStateDof \modelOpDof_{n+1|n})  \projCovDof_n + \diffParamDof \modelOpDof_{n+1|n}
   	\end{cases}
   \end{align}
	\item Correction:
   \begin{align}
   	\begin{cases}	
   		\redCovDof_{n+1} &= \redCovDof_{n} + (\projStateCovDof_{n+1})^\intercal (\diffStateDof  \obsOpDof_{n+1})^\intercal \obsNoiseNormDof_{n+1} (\diffStateDof \obsOpDof_{n+1}) \projStateCovDof_{n+1} \\
		\hat{\stateDof}_{n+1}^+  &= \hat{\stateDof}_{n+1}^- + \projStateCovDof_{n+1} (\hat{\paramDof}_{n+1}^+  - \hat{\paramDof}_{n+1}^-) \\
		\hat{\paramDof}_{n+1}^+  &= \hat{\paramDof}_{n+1}^- \\
		& \hspace{0.1 cm} +  \redCovDof_{n+1}^{-1}(\projStateCovDof_{n+1})^\intercal(\diffStateDof  \obsOpDof_{n+1})^\intercal\obsNoiseNormDof_{n+1}(\observDof_{n+1} - \obsOpDof_{n+1}(\hat{\stateDof}_{n+1}^-))
   	\end{cases}
   \end{align}
\end{enumerate}
\end{subequations}
starting from the same initial conditions.

We have presented the Reduced Order EKF and have recalled its convergence. Now a legitimate question is to know if a reduced order strategy can be applied to other types of approximation of the optimal filtering approach. The answer was given for the UKF filter in \cite{MoireauReducedorder} where a reduced order version (RoUKF) was derived. Moreover it was also applied  to a specific case where  the initial uncertainty is reduced to the parametric space. In the general context --~\emph{i.e.} without particularizing the state and parameter dependency --~the RoUKF can be formulated as follow.

Given an adequate sampling rule, we store the corresponding weights $(\alpha_i)$ in the diagonal matrix $\weightUKFmat_\alpha$ and precompute specific unitary simplex sigma-points $(I^{i})_{1\leq i \leq \Nsigma}$ (i.e.~with zero mean and unit covariance) with $\Nsigma = \Nparam+1$ since the reduced space corresponds here to the parametric space. The algorithm consists of the following three steps computed recursively:
\begin{subequations}\label{eq:ROUKF}	
	\begin{enumerate}
		\item Sampling: 
	\begin{align}
		\begin{cases}
			\redCovDofSqrt_n &= \sqrt{(\redCovDof_n)^{-1}} \\[0.1cm]
			\hat{\stateDof}^{[i]+}_{n} &= \hat{\stateDof}_{n}^+ + \projStateCovDof_n \cdot \redCovDofSqrt_n^\intercal \cdot I^{i},\quad 1\leq i \leq \Nsigma \\[0.1cm]
						\hat{\paramDof}^{[i]+}_{n} &= \hat{\paramDof}_{n}^+ + \projParamCovDof_n \cdot \redCovDofSqrt_n^\intercal \cdot I^{i},\quad 1\leq i \leq \Nsigma \\[0.1cm]
		\end{cases}
	\end{align}
	\item Prediction:
   \begin{align} 
   	\begin{cases}			
   		\hat{\stateDof}^{[i]-}_{n+1} &= \modelOpDof_{n+1|n}(\hat{\stateDof}^{[i]+}_{n},\hat{\paramDof}^{[i]+}_{n}),\quad 1\leq i \leq \Nsigma \\[0.1cm] 
   		\hat{\paramDof}^{[i]-}_{n+1} &= \hat{\paramDof}^{[i]+}_{n},\quad 1\leq i \leq \Nsigma \\[0.1cm] 
   	   	\hat{\stateDof}^-_{n+1} &= \sum_{i=1}^{\Nsigma} \alpha_i \hat{\stateDof}_{n+1}^{[i]-}\\[0.1cm]
   	   	\hat{\paramDof}^-_{n+1} &= \sum_{i=1}^{\Nsigma} \alpha_i \hat{\paramDof}_{n+1}^{[i]-}\\[0.1cm]
   	\end{cases}
   \end{align}
	\item Correction:
   \begin{align}
   	\begin{cases}	
   		\projStateCovDof_{n+1} &= [\hat{\stateDof}^{[*]-}_{n+1}]\weightUKFmat_\alpha [I^{[*]}]^\intercal \\[0.1cm]		
	\projParamCovDof_{n+1} &= [\hat{\paramDof}^{[*]-}_{n+1}]\weightUKFmat_\alpha [I^{[*]}]^\intercal \\[0.1cm]		   		
   		\observDof_{n+1}^{[i]-} &= \obsOpDof_{n+1}(\hat{\stateDof}^{[i]-}_{n+1}) \\[0.1cm]
		\observDof_{n+1}^{-} &= \sum_{i=1}^{\Nsigma} \alpha_i \observDof_{n+1}^{[i]-} \\[0.1cm]
   		\covObsUKF_{n+1} &=  [\observDof^{[*]-}_{n+1}]\weightUKFmat_\alpha [I^{[*]}]^\intercal \\[0.1cm]
   		\redCovDof_{n+1} &=  \1 + \covObsUKF_{n+1}^\intercal  \obsNoiseNormDof_{n+1} \covObsUKF_{n+1} \\[0.1cm]
   		\hat{\stateDof}_{n+1}^{+} &= \hat{\stateDof}_{n+1}^{-} + \projStateCovDof_{n+1} \redCovDof_{n+1}^{-1} \covObsUKF_{n+1}^\intercal \obsNoiseNormDof_{n+1}  (\observDof_{n+1} - \observDof_{n+1}^{-})\\[0.1cm]
\hat{\paramDof}_{n+1}^{+} &= \hat{\paramDof}_{n+1}^{-} + \projParamCovDof_{n+1} \redCovDof_{n+1}^{-1}  \covObsUKF_{n+1}^\intercal \obsNoiseNormDof_{n+1}  (\observDof_{n+1} - \observDof_{n+1}^{-})   		
   	\end{cases}
   \end{align}
\end{enumerate}
\end{subequations}
where $[I^{[*]}]$ is the matrix concatenating the $(I^{i})$ vectors side by side, and similarly for other vectors.

The analysis of this estimator can also be found in \cite{MoireauReducedorder} and relies again on the fact that the linearized error satisfies exactly the same dynamics as the linearized error of the RoEKF. Therefore the demonstration by linearization is directly obtained. The combined practical simplicity and efficiency of such parameter estimator in comparison to other adaptive observers have made this approach popular for real case parameter estimation problems  \cite{bertoglio-moireau-gerbeau-11,chabiniok-moireau-ea-11,marchesseau:hal-00819806,Xi-ea-10}.

\section*{References}

\bibliographystyle{plain}

\end{document}